\documentclass[notitlepage,leqno,10pt]{article}

\usepackage{latexsym}
\usepackage{amsmath}
\usepackage{amsfonts}
\usepackage{amssymb}
\renewcommand{\a }{\alpha }
\renewcommand{\b }{\beta }
\renewcommand{\d}{\delta }
\newcommand{\D }{\Delta }
\newcommand{\Di }{\mathcal{D}^{1,2}(\R^n_+)}
\newcommand{\Dit }{\mathcal{D}^{1,2}(\R^3_+)}
\newcommand{\e }{\varepsilon }
\newcommand{\g }{\gamma}
\renewcommand{\i }{\iota}
\newcommand{\G }{\Gamma }
\renewcommand{\l }{\lambda }

\newcommand{\n }{\nabla }
\newcommand{\var }{\varphi }
\newcommand{\rh }{\rho }
\newcommand{\vrh }{\varrho }
\newcommand{\s }{\sigma }
\newcommand{\Sig }{\Sigma}
\renewcommand{\t }{\tau }
\renewcommand{\th }{\theta }
\renewcommand{\o }{\omega }
\renewcommand{\O }{\Omega }

\newcommand{\ov}{\overline}

\newcommand{\be}{\begin{equation}}
\newcommand{\ee}{\end{equation}}
\newenvironment{pf}{\noindent{\sc Proof}.\enspace}{\rule{2mm}{2mm}\medskip}
\newenvironment{pfn}{\noindent{\sc Proof}}{\rule{2mm}{2mm}\medskip}

\newcommand{\R}{\mathbb{R}}

\newcommand{\N}{\mathbb{N}}

\author{Zindine DJADLI$^{\rm a}$ \and Andrea MALCHIODI$^{\rm b}$ \and Mohameden OULD
AHMEDOU$^{\rm c}$}

\date{}

\title{PRESCRIBING SCALAR AND BOUNDARY MEAN \\ CURVATURE ON THE
THREE DIMENSIONAL \\ HALF SPHERE}

\begin{document}

\newtheorem{lem}{Lemma}[section]
\newtheorem{pro}[lem]{Proposition}
\newtheorem{thm}[lem]{Theorem}
\newtheorem{rem}[lem]{Remark}
\newtheorem{cor}[lem]{Corollary}
\newtheorem{df}[lem]{Definition}

\maketitle

\begin{center}

$^{\rm a}${\small
Universit\'e de Cergy-Pontoise, D\'epartement de Math\'ematique,
\\ Site de Saint-Martin, 2 Avenue Adolphe Chauvin, F 95302
Cergy-Pontoise Cedex, France}

$^{\rm b}${\small School of Mathematics, Institut for Advanced Studies,
\\ 1 Einstein Drive, 088540
Princeton, NJ, USA}

$^{\rm c}${\small Mathematisches Institut der Universit\"at Bonn, \\
 Beringstrasse 4,D-53115 Bonn, Germany}

\end{center}

\footnotetext[1]{E-mail addresses: Zindine.Djadli@math.u-cergy.fr
(Z. Djadli), malchiod@ias.edu (A. Malchiodi),
ahmedou@ahmedou@math.uni-bonn.de (M. Ould Ahmedou).}

\noindent {\sc abstract}. - We consider the problem of prescribing
the scalar curvature and the boundary mean curvature of the
standard half three sphere, by deforming conformally its standard
metric. Using blow up analysis techniques and minimax arguments,
we prove some existence and compactness results.

%

\section{Introduction}\label{s:i}
In this paper we study some equation arising in differential
geometry, when the metric of a Riemannian manifold is conformally
deformed. Precisely, is given a manifold with boundary $(M,g)$ of
dimension $n \geq 3$; transforming the metric $g$ into $g' =
v^{\frac{4}{n-2}}g$, where $v$ is a smooth positive function, the
scalar curvatures $R_g, R_{g'}$ and the mean curvatures of the
boundary $h_g, h_{g'}$, with respect to $g$ and $g'$ respectively,
are related by the formulas
\begin{equation}\label{eq:in}
\begin{cases}
-4\frac{n-1}{n-2} \Delta_g v + R_g v = R_{g'} v^{\frac{n+2}{n-2}}, &
\mbox{ in } M; \\
\frac{2}{n-2} \frac{\partial v}{\partial \nu} + h_g v = h_{g'}
v^{\frac{n}{n-2}}, & \mbox{ on } \partial M,
\end{cases}
\end{equation}
see e.g. \cite{aul}. In the above equation, $\nu$ denotes the outward
unit vector perpendicular to $\partial M$, with respect to the metric
$g$.

A problem arises naturally when looking at equation \eqref{eq:in}:
assigned two functions $K : M \to \R$ and $H : \partial M \to \R$,
does exists a metric $g'$ conformally equivalent to $g$ such that
$R_{g'} \equiv K$ and $h_{g'} \equiv H$? From equation
\eqref{eq:in}, the problem is equivalent to finding a smooth
positive solution $v$ of the equation
\begin{equation}\label{eq:in2}
\begin{cases}
-4\frac{n-1}{n-2} \Delta_g v + R_g v = K v^{\frac{n+2}{n-2}}, &
\mbox{ in } M; \\
\frac{2}{n-2} \frac{\partial v}{\partial \nu} + h_g v = H
v^{\frac{n}{n-2}}, & \mbox{ on } \partial M.
\end{cases}
\end{equation}
The requirement about the positivity of $v$ is necessary for the
metric $g'$ to be Riemannian. For the two-dimensional case, there
are analogous equations involving exponential nonlinearities.

We are mainly interested in the so-called {\em positive case}, see
\cite{E1}, when the quadratic part of the Euler functional associated
to \eqref{eq:in2} is positive definite.

A first criterion for existence of solutions of \eqref{eq:in2},
and also a proof of regularity, was given by P. Cherrier,
\cite{ch}. He proved that if the energy of some test function is
below an explicit threshold, then problem \eqref{eq:in2} admits a
solution as a mountain pass critical point. Using this criterion
J. Escobar obtained some existence results in the interesting
particular case of constant $K$ and $H$, see \cite{E1}, \cite{E2}.
The proof relies on some extension of the Positive Mass Theorem by
R. Schoen and S.T. Yau \cite{sy} to the case of manifolds with
boundary. He showed that almost every compact manifold with
boundary can be conformally deformed so that its scalar curvature
is $1$ and the boundary is minimal, i.e. the mean curvature is
zero. He also gives some results when $H$ is a constant close to
zero. More recently Z.C. Han and Y.Y. Li, see \cite{HL1},
\cite{HL2}, extended most of the results of J. Escobar to the case
in which $K \equiv 1$ and $H$ is any constant. They also prove a
compactness results in the locally conformally flat manifolds with
umbilic boundary.

\

\noindent We consider here the case of the standard half sphere
$S^n_+ = \{ x \in \R^{n+1} \, : \, \|x\| = 1, x_{n+1} > 0\}$
endowed with its standard metric $g_0$, and in particular the case
$n = 3$; the functions $K$ and $H$ are now non constant and $K$
will always assumed to be positive. We are thus reduced to find
positive solutions $v$ of the problem
\begin{equation}\label{eq:sp}
\begin{cases}
- 4\frac{n-1}{n-2} \, \D v + n (n-1) \, v =
K(x) \, v^{\frac{n+2}{n-2}}, & \mbox{ in } S^n_+; \\
\frac{2}{n-2}\frac{\partial v}{\partial \nu} = H(x') \,
v^{\frac{n}{n-2}}, & \mbox{ on } \partial S^n_+.
\end{cases}
\end{equation}

Problem \eqref{eq:sp} is in some sense related to the well-known {\em
Scalar Curvature Problem} on $S^n$
\begin{equation}\label{eq:scsn}
- 4\frac{n-1}{n-2} \, \D v + n (n-1) \, v = K(x) \,
v^{\frac{n+2}{n-2}}, \qquad \mbox{ in } S^n,
\end{equation}
to which much work has been devoted, see \cite{agp}, \cite{bf},
\cite{bc}, \cite{bcch}, \cite{cl}, \cite{[ACGPY]}, \cite{[ACPY1]},
\cite{[ACPY2]}, \cite{[ACPY3]}, \cite{hv}, \cite{yy1}, \cite{yy2},
\cite{sz} and references therein. As for \eqref{eq:scsn}, also for
problem \eqref{eq:sp} there are topological obstructions for
existence of solutions, based on Kazdan-Warner type conditions,
see \cite{bp} and also the proof of Proposition \ref{p:cto}. Hence
it is not expectable to solve problem \eqref{eq:sp} for all the
functions $K$ and $H$, and it is natural to impose some conditions
on them.

We would like to point out the following features of the scalar
curvature problem on $S^n$ in lower dimensions, referring to the
above-mentioned papers. For $n = 2$, non-converging Palais Smale
sequences are characterized by the presence of just one bubble.
Under generic assumptions on $K$, it turns out that when $n = 3,
4$, solutions of \eqref{eq:scsn}, or of some subcritical
approximation, possess only {\em isolated simple blow ups}, see
Section \ref{s:bp} for the definition. When $n = 3$ there is
indeed just one blow up point, while for $n = 4$ blow ups may as
well occur at more points.

We now discuss problem \eqref{eq:sp}. For $n = 2$ (about the
corresponding equation with exponential nonlinearities), P.L. Li
and J.Q. Liu proved in \cite{ll} that compactness is lost along
one bubble only, as in the case of the problem on the sphere. The
only difference is that blow up can only occur at the boundary of
$S^2_+$. For the case $H \equiv 0$ and positive $K$, they prove
existence results which are in some sense reminiscent of those of
\cite{cl}, see also \cite{[ACPY1]}, \cite{[ACPY2]}.

In \cite{yn} Y.Y. Li considered the case of $n = 3$ and $H \equiv
0$. Under generic assumptions on $K$, he proved that blow ups are
isolated simple and at only the boundary. He also stated that, as
for \eqref{eq:scsn} when $n = 3$, blow ups may occur at most one
point. Actually the last statement is not true, although the main
features of the blow-up behavior at the boundary are analyzed in
\cite{yn}. Using the ingredients of \cite{HL1}, \cite{yn} and
\cite{LZ1} we correct this here and we prove that (also for non
constant $H$) blow ups, which are always isolated simple and on
the boundary, can be multiple, see Section \ref{s:6}. Hence the
situation could be considered similar to that of \eqref{eq:scsn}
for $n = 4$. This fact, in the case of $H \equiv 0$, could be
roughly explained as follows. Reflecting both $K$ and $v$ evenly
to the whole $S^3$, one could study symmetric solutions of
\eqref{eq:scsn}. The blow up analysis for the three dimensional
case strongly relies on the differentiability of $K$ at blow up
points. This implies that blow ups of symmetric functions outside
$\partial S^3_+$, which are multiple, are ruled out. This argument
does not apply when blow up points are on $\partial S^3_+$, since
the symmetric extension of $K$ is not regular there. See Remark
\ref{r:mu} for a more quantitative explanation of this fact.

For the case of any $n$, some results are proved in \cite{alm}
when $K$ and $H$ are close to some constants; here we are
extending some of those results for $n = 3$ without the {\em close
to constant} conditions, see also Remark \ref{r:n}. In the paper
\cite{cxy} the case $K \equiv 0$ and $H$ close to a positive
constant is considered, for $n \geq 3$. In the forthcoming paper
\cite{dmo} we will extend some of those results to the
non-perturbative case.

Our first result is the following.

\begin{thm}\label{t:mp}
Assume $n = 3$, let $K : \ov{S^3_+} \to \R$ be a positive $C^1$
function and let $H : \partial S^3_+ \to \R$ be of class $C^2$.
Let $\varphi : \partial S^3_+ \to \R$ be the function defined by
\begin{equation}\label{eq:vph}
\varphi(x') = 4 \pi \, \sqrt{\frac{6}{K(x')}} \, \left( \frac{\pi}{2}
- \arctan\left( H(x') \, \sqrt{\frac{6}{K(x')}} \right) \right),
\qquad x' \in \partial S^3_+.
\end{equation}
Suppose that for some point $q \in \partial S^3_+$ the following
condition holds
\begin{equation}\label{eq:p-max}
\varphi(q) = \min_{\partial S^3_+} \varphi; \qquad \qquad
\frac{\partial K}{\partial \nu}(q) < 0.
\end{equation}
Then there exists a positive solution of problem \eqref{eq:sp}.
\end{thm}
The proof of Theorem \ref{t:mp} relies on the study of the following
subcritical approximation of equation \eqref{eq:sp}
\begin{equation}\label{eq:spa}
\begin{cases}
- 4\frac{n-1}{n-2} \, \D u + n (n-1) \, u =
K \, u^{p}, & \mbox{ in } S^n_+; \\
\frac{2}{n-2}\frac{\partial u}{\partial \nu} = H \,
u^{\frac{p+1}{2}}, & \mbox{ on } \partial S^n_+.
\end{cases}
\end{equation}
Here the exponent $p$ is converging to $\frac{n+2}{n-2}$ from
below. As mentioned before, for $n = 3$ blow ups of equation
\eqref{eq:spa} can occur only at the boundary of $S^3_+$.
Nevertheless, if $v_p$ denotes a mountain pass solution of
\eqref{eq:spa} for $p < \frac{n+2}{n-2}$, condition
\eqref{eq:p-max} implies that $\{v_p\}_p$ is uniformly bounded for
$p \to \frac{n+2}{n-2}$, and hence converges to a solution of
\eqref{eq:sp}. The function $\varphi(x')$ represents the blow up
energy at a point $x' \in
\partial S^3_+$ and plays a crucial role in the blow up analysis.
Indeed, see Section \ref{s:6}, blow ups can occur only at critical
points of $\varphi$. We note that when $H$ is a constant function,
critical points of $\varphi$ coincide with critical points of
$K|_{\partial S^3_+}$, see also \cite{yn}.

Under generic assumptions on $K$ and $H$, ($(K,H) \in \mathcal{A}$
in the notation below), it is possible to give a complete
description of the behavior of general solutions of \eqref{eq:spa}
when $p$ converges to $\frac{n+2}{n-2}$, and to deduce existence
and compactness results for equation \eqref{eq:sp}. We point out
that, in order to this, we use crucially the classification result
in \cite{LZ1} and the blow-up analysis in \cite{HL1}, \cite{yy2}.
The blow up analysis provides necessary conditions on these
solutions, while the Implicit Function Theorem gives sufficient
conditions for existence of solutions highly concentrating at some
points of $\partial S^3_+$. In this way one can compute the total
Leray Schauder degree of the solutions of \eqref{eq:sp} in the
space $C^{2,\a}(\ov{S^3_+})$, for some $\a \in (0,1)$. Such a
method has been used in \cite{sz} and \cite{yy2} for problem
\eqref{eq:scsn} in dimensions $3$ and $4$ respectively.

To state our next result we need to introduce some notation, which
considerably simplifies in the case $H \equiv 0$, see Remark
\ref{r:h0}. Given $K \in C^2(\ov{S^3_+})$ and $H \in C^2(\partial
S^3_+)$, let $\varphi \in C^2(\partial S^3_+)$ be defined by formula
\eqref{eq:vph}, and set
$$
\mathcal{F} = \left\{ q \in \partial S^3_+ \, : \, \n \varphi(q) = 0
\right\}; \qquad \mathcal{F}^{+ (-)} = \left\{ q \in \partial S^3_+
\, : \, \n' \varphi(q) = 0, \frac{\partial K}{\partial \nu}(q) > 0 \,
(< 0) \right\};
$$
$$
\mathcal{M}_{K,H} = \left\{ v \in C^{2,\a}(\ov{S^3_+}) \, : \, v
\mbox{ satisfies} \eqref{eq:sp} \right\}.
$$
Here $\n'$ denotes the gradient of functions defined on $\partial
S^n_+$.

For $q \in \partial S^3_+$, let $\pi_{q} : S^3_+ \to \R^3_+$ denote
the stereographic projection with pole $- q$. In
$\pi_{q}$-stereographic coordinates, we define the function $G_{q} :
S^3_+ \to \R$ by
\begin{equation}\label{eq:JP0}
G_{q}(x) = \left( \frac{1 + |x|^2}{2} \right)^{\frac{1}{2}}
\frac{1}{|x|}, \qquad x \in \R^3_+.
\end{equation}
The function $G_{q}$ is the Green's function for the conformal
laplacian $- 8 \D + 6$ on $S^3$ with pole $q$. Define also $\psi :
\partial S^3_+ \to \R$ by
\begin{equation}\label{eq:psi}
\psi(x') = 1 + H(x') \sqrt{\frac{6}{K(x')}} \left( \arctan \left(
H(x') \sqrt{\frac{6}{K(x')}} \right) - \frac{\pi}{2} \right) = 1 -
\frac{H(x') \varphi(x')}{4 \pi}.
\end{equation}
To each $\{ q^1, \dots, q^N \} \subseteq \mathcal{F} \setminus
\mathcal{F}^-$, $N \geq 1$, we associate an $N \times N$ symmetric
matrix   $M = M(q^1, \dots, q^N)$ defined by
\begin{equation}\label{eq:M}
\begin{cases}
M_{jj} = \frac{\partial K}{\partial \nu}(q^j) \,
\frac{\psi(q^j)}{K(q^j)^{\frac{3}{2}}}, &
j \in \{1, \dots, N\} \\
M_{lj} = - 4 \sqrt{2} \, \frac{G_{q^l}(q^j)}{K(q^l)^{\frac{1}{4}}
K(q^j)^{\frac{1}{4}}}, & l, j \in \{1, \dots, N\}, \quad l \neq j.
\end{cases}
\end{equation}
Let $\rho = \rho(q^1, \dots, q^N)$ denote the least eigenvalue of
$M$.
It has been first pointed out by A. Bahri, \cite{bab}, see also
\cite{bc}, that when the interaction between different bubbles is
of the same order as the {\em self interaction}, the function
$\rho$ for a matrix as in \eqref{eq:M} plays a fundamental role in
the theory of the critical points at infinity. For problem
\eqref{eq:sp}, such kind of phenomenon appears when $n = 3$.

Define the set $\mathcal{A}$ to be
\begin{eqnarray*}
\mathcal{A} & = & \{ (K, H) \in C^2(\ov{S^3_+}) \times C^2(\partial
S^3_+) \, : \, K > 0, \varphi \mbox{ is a Morse function
on $\partial S^3_+$, } \\
& & \frac{\partial K}{\partial \nu} \neq 0 \hbox{ on }
\mathcal{F}, \mbox{ and } \rho = \rho(q^1, \dots, q^N) \neq 0,
\forall \, q^1, \dots, q^N \in \mathcal{F} \}.
\end{eqnarray*}
Let us observe that the condition $(K, H) \in \mathcal{A}$ is
generic. We introduce an integer valued continuous function Index $:
\mathcal{A} \to \N$ by the following formula
$$
\mbox{Index}(K,H) = -1 + \sum_{j=1}^\ell
\sum_{\stackrel{\rho(q^{i_1}, \dots q^{i_j}) > 0,}{1 \leq i_1 < i_2
\leq \dots \leq i_j \leq \ell}} (-1)^{j + \sum_{l=1}^j (2 -
m(\varphi, q^j))},
$$
where $m(\varphi,q^{i_l})$ denotes the Morse index of $\varphi$ at
$q^{i_l}$, and $\ell = card \, | \mathcal{F^+}|$. Now we are able to
state our next result, about existence and compactness of solutions
of \eqref{eq:sp}.

\begin{thm}\label{t:61}
Let $n = 3$  and suppose $(K,H) \in \mathcal{A}$. Then for all $
\a \in (0,1)$, there exists some constant $R$ depending only on
$\min_{\ov{S^3_+}} K$, $\|K\|_{C^1(\ov{S^3_+})}$,
$\|H\|_{C^2(S^3_+)}$, $\min \{ |\rho (q^1, \dots, q^N)| \, : \,
q^1, \dots, q^N \in \mathcal{F}, N \geq 2 \}$ and $\a$ such that
\begin{equation}\label{eq:cc6}
\frac{1}{R} \leq v \leq R, \qquad \|v\|_{C^{2,\a}(\ov{S^3_+})} \leq
R,
\end{equation}
for all positive solutions $v$ of equation \eqref{eq:sp}.
Moreover \eqref{eq:sp} possesses a solution provided {\em
Index}$(K,H) \neq 0$.
\end{thm}

Since the situation here resembles that of $S^4$ for a Morse
function $K$, our Theorem \ref{t:61} can be considered as a
counterpart of the results in \cite{bcch} and \cite{yy2} for
manifolds with boundary. Notice that only the least eigenvalue of
$M(q^1, \dots, q^N)$ plays a role in counting the total degree of
solutions of \eqref{eq:sp} and in the compactness. For instance,
considering a continuous family of functions $(K_t, H_t)$, the
total degree changes when the least eigenvalue of $M_t(q^1, \dots,
q^N)$ crosses zero, while it remains unchanged when other
eigenvalues cross zero.

\begin{rem}\label{r:h0}
(a) When the function $H$ is identically equal to zero, the
functions $\varphi(x')$ and $\psi(x')$ assume the simpler form
$$
\varphi(x') = 2 \pi^2 \sqrt{\frac{6}{K(x')}}; \qquad \qquad \psi(x')
= 1.
$$
In particular minima of $\varphi$ coincide with maxima of $K$
restricted to the boundary and viceversa.  (b) We note that
equation \eqref{eq:sp} for $n = 3$ is invariant under the
rescaling
$$
K \to \g K; \qquad H \to \g^{\frac{1}{2}} H; \qquad u \to
\g^{-\frac{1}{4}} u,
$$
where $\g$ is any positive number. The hypotheses involving $K$ and
$H$ in Theorems \ref{t:mp} and \ref{t:61} are both invariant under
such a rescaling.

(c) While Theorem \ref{t:61} is related to some known results for
equation \eqref{eq:scsn}, Theorem \ref{t:mp} has no counterpart in
the problem on the whole sphere. The existence argument is
strictly related to the presence of the boundary.
\end{rem}

\noindent The authors have been recently informed about some
related results obtained in \cite{ga}.

\

\noindent The paper is organized as follows. In Section 2 we
collect some useful technical tools, while in Section 3 we compute
the blow up energies, depending on the values of $K$ and $H$ at
the blow up point. In section 4 we recall some known facts about
blow up analysis of equations \eqref{eq:sp} and \eqref{eq:scsn},
and in section 5 we specialize to the case of boundary blow ups.
Then in Section 6 we prove that blow ups are isolated simple, see
Definition \ref{d:isbu}, and occur only at the boundary of
$\partial S^3_+$. Finally, Sections 7 and 8 are devoted to the
proofs of Theorems \ref{t:mp} and \ref{t:61} respectively.

\begin{center}
{\bf Acknowledgments}
\end{center}

\noindent This work was first began while Z. D. enjoyed the
hospitality of SISSA at Trieste. He would like to thank A.
Ambrosetti for the invitation. Part of it was accomplished when Z.
D. was visiting the Mathematical Sciences Research Institute
(MSRI) at Berkeley. His research at MSRI is supported in part by
NSF grant DMS-9701755. A. M. is supported by a Fulbright
fellowship for the academic year 2000-2001 and by M.U.R.S.T.,
under the project {\em Variational Methods and Nonlinear
Differential Equations}. M.O.A.'s research has been  supported by a
S.I.S.S.A. postdoctoral fellowship.

\section{Some preliminaries}\label{s:p}

We will use the notation $x$ for variables belonging to the half
sphere $S^n_+$, or to the half space $\R^n_+$, which is defined by
$\R^n_+ := \{ x \in \R^n \, : \, x_n  > 0 \}$; variables in both the
boundaries will be denoted in general by $x'$.

Solutions of problem \eqref{eq:sp} are critical points of the Euler
functional $J_{K,H} : H^1(S^n_+) \to \R$ defined in the following way
\begin{eqnarray}\label{eq:jkh}
J_{K,H}(v) = \frac{1}{2} \int_{S^n_+} \left( \frac{4(n-1)}{n-2} |\n
v|^2 + n(n-1) \, v^2 \right) - \frac{1}{2^*} \int_{S^n_+} K(x) \,
|v|^{2^*} - (n-2) \int_{\partial S^n_+} H(x') \,
|v|^{2\frac{n-1}{n-2}}, \nonumber \\ v \in H^1(S^n_+),
\end{eqnarray}
where $2^* = \frac{2n}{n-2}$. It will be convenient to perform some
stereographic projection in order to reduce the above problem to
$\R^n_+$. Let $\Di$ denote the completion of
$C_c^\infty\left(\ov{\R^n_+}\right)$ with respect to the Dirichlet
norm. The stereographic projection $\pi_q$ through a point $q \in
\partial S^n_+$ induces an isometry $\iota : H^1(S^n_+) \to \Di$
according to the following formula
\begin{equation}\label{eq:iota}
(\iota \, v)(x) = \left( \frac{2}{1 + |x|^2} \right)^{\frac{n-2}{2}}
\, v(\pi_q^{-1}(x)), \qquad \qquad v \in H^1(S^n_+), \quad x \in
\R^n_+.
\end{equation}
In particular, one can check that the following relations hold true,
for every $v \in H^1(S^n_+)$
$$
\int_{S^n_+} \left(\frac{4(n-1)}{n-2} |\n v|^2 + n(n-1) \, v^2
\right) = \int_{\R^n_+}\frac{4(n-1)}{n-2}  |\n (\iota \, v)|^2;
$$
$$
\int_{S^n_+} |v|^{2^*} = \int_{\R^n_+} |\i \, v|^{2^*}; \qquad
\int_{\partial S^n_+} |v|^{2\frac{n-1}{n-2}} = \int_{\partial \R^n_+}
|\i \, v|^{2\frac{n-1}{n-2}}.
$$
By means of these equations, the functional $J_{K,H}$ transforms into
$I_{K,H} : \Di \to \R$ given by
$$
I_{K,H}(u) = \frac{1}{2} \int_{\R^n_+} \frac{4(n-1)}{n-2} |\n u|^2 -
\frac{1}{2^*} \int_{\R^n_+} K(x) \, |u|^{2^*} - (n-2) \int_{\partial
\R^n_+} H(x') \, |u|^{2\frac{n-1}{n-2}},
$$
meaning that $J_{K,H}(v) = I_{K,H}(\iota \, v)$ for every $v \in
H^1(S^n_+)$. Here we are identifying the functions $K$ and $H$ and
their compositions with the stereographic projection $\pi_q$. This
fact will be assumed as understood in the sequel.

Critical points of the functional $I_{K,H}$ are solutions of the
following problem
\begin{equation}\label{eq:rp}
\begin{cases}
- 4\frac{n-1}{n-2} \, \D u = K(x) \, u^{\frac{n+2}{n-2}}, &
\mbox{ in } \R^n_+; \\
- \frac{2}{n-2}\frac{\partial u}{\partial x_n} = H(x') \,
u^{\frac{n}{n-2}}, & \mbox{ on } \partial \R^n_+.
\end{cases}
\end{equation}
Using similar arguments, one finds that the counterpart in $\R^n_+$
of equation \eqref{eq:spa} is given by
\begin{equation}\label{eq:rpp}
\begin{cases}
- 4\frac{n-1}{n-2} \, \D u = W(x)^\t \, K(x) \, u^{p},
& \mbox{ in } \R^n_+; \\
- \frac{2}{n-2}\frac{\partial u}{\partial x_n} = W(x')^{\frac{\t}{2}}
H(x') \, u^{\frac{p+1}{2}}, & \mbox{ on } \partial \R^n_+,
\end{cases}
\end{equation}
where $W(x) = \left( \frac{2}{1+|x|^2} \right)^{\frac{n-2}{2}}$, and
where $\t = \frac{n+2}{n-2} - p$. The terms $W(x)^\t$ and
$W(x')^{\frac{\t}{2}}$ in the equation above are corrections due to
the non conformality of equation \eqref{eq:spa} when $p \neq
\frac{n+2}{n-2}$.

As a typical feature of non compact variational problems like
\eqref{eq:rp}, it is fundamental to analyze the associated {\em
problems at infinity}. Solutions of such problems describe the
asymptotic profile of non-converging Palais Smale sequences. In the
specific case of \eqref{eq:rp}, these problems at infinity are of two
kinds, namely
\begin{equation}\label{eq:pirn}
- 4\frac{n-1}{n-2} \, \D u = K(\ov{x}) \, u^{\frac{n+2}{n-2}}, \qquad
\hbox{ in } \R^n,
\end{equation}
for some fixed $\ov{x} \in \ov{\R^n_+}$, and
\begin{equation}\label{eq:pirn+}
\begin{cases}
- 4\frac{n-1}{n-2} \, \D u = K(\ov{x}') \, u^{\frac{n+2}{n-2}}, &
\mbox{ in } \R^n_+; \\
- \frac{2}{n-2}\frac{\partial u}{\partial x_n} = H(\ov{x}') \,
u^{\frac{n}{n-2}}, & \mbox{ on } \partial \R^n_+,
\end{cases}
\end{equation}
for some $\ov{x}' \in \partial \R^n_+$. Roughly, problem
\eqref{eq:pirn} corresponds to the case in which the functions are
mostly concentrated in the interior of $\R^n_+$, while problem
\eqref{eq:pirn+} corresponds to the case in which the functions are
concentrated near the boundary. We note that solutions of problem
\eqref{eq:pirn+} are critical points of the functional
$$
I_{\ov{K},\ov{H}}(u) = \frac{1}{2} \int_{\R^n_+} \frac{4(n-1)}{n-2}
|\n u|^2 - \frac{1}{2^*} \ov{K} \int_{\R^n_+} |u|^{2^*} - (n-2)
\ov{H} \int_{\partial \R^n_+} |u|^{2\frac{n-1}{n-2}},
$$
where $\ov{K} = K(\ov{x}')$ and $\ov{H} = H(\ov{x}')$.

Positive solutions of problems \eqref{eq:pirn} and
\eqref{eq:pirn+} have been completely classified in \cite{cgs} and
\cite{LZ1}, see also \cite{ei}; we recall the results in the
following Lemma.

\begin{lem}\label{l:cs}
The positive solutions $u$ of problems \eqref{eq:pirn} and
\eqref{eq:pirn+} are, modulo translations in $\R^n$ and $\R^{n-1}$,
respectively of the form
\begin{equation}\label{eq:ulrn}
U_\l(x) = \left( \frac{\l}{1 + \l^2 k |x|^2}
\right)^{\frac{n-2}{2}}; \qquad \l > 0,
\end{equation}
where $k = \frac{K(\ov{x})}{4n(n-1)}$, and
\begin{equation}\label{eq:ulrn+}
\ov{U}_\l(x) = \left( \frac{\l}{1 + \l^2 k \left( |x'|^2 +
(x_n + t_\l)^2 \right)} \right)^{\frac{n-2}{2}}; \qquad \l > 0,
\end{equation}
where $k = \frac{K(\ov{x}')}{4n(n-1)}$ and where $t_\l$ is given by
$2 \, k \, t_\l \, \l = H(\ov{x}')$.
\end{lem}

\noindent It will be convenient to consider the expression
$\int_{\R^n_+} x_n \ov{U}_1^{\frac{2n}{n-2}}$ with $n = 3$. Using the
formula
$$
\int_0^\infty \frac{r^\a}{(1+r^2)^\b} =
\frac{\G\left(\frac{\a+1}{2}\right) \G\left(\b -
\frac{\a+1}{2}\right)}{2 \G\left(\b\right)},
$$
and integrating first in the variable $x'$ and then in the variable
$x_3$, one obtains
\begin{equation}\label{eq:intxnu}
    \int_{\R^3_+} x_3 \ov{U}_1^{6} = \frac{144 \pi}{K(\ov{x}')^2} \left[ 1 +
    H(\ov{x}') \sqrt{\frac{6}{K(\ov{x}')}} \left( \arctan
    \left(H(\ov{x}')\sqrt{\frac{6}{K(\ov{x}')}} - \frac{\pi}{2} \right)
    \right)\right] = \frac{144 \pi}{K^2(\ov{x}')} \, \psi(\ov{x}').
\end{equation}

\noindent The projection to $\R^n_+$ is not the only
transformation we will perform. In the next section we will use a
conformal transformation of $S^n_+$ onto some suitable spherical
cap $\Sig_\th$. In a similar way as before, this transformation
induces an isometry between $H^1(S^n)$ and $H^1(\Sig_\th)$. We
will not write the explicit formulas for this transformation, as
in \eqref{eq:iota}.

\section{Study of the blow up energies}\label{s:bue}

In this section we compute the energies of the solutions of problems
\eqref{eq:pirn} and \eqref{eq:pirn+} (i.e. of the functions given in
\eqref{eq:ulrn} and \eqref{eq:ulrn+}), highlighting the dependence on
the values of $K$ (and $H$) at $\ov{x}$ (resp. at $\ov{x}'$). It is
well known that these energies are strongly related to those of non
converging Palais Smale sequences of the functionals $J_{K,H}$ and
$I_{K,H}$.

In order to simplify the computations it is convenient, using a
suitable stereographic transformations, to reduce problems
\eqref{eq:pirn} and \eqref{eq:pirn+} to $S^n$ and to some
spherical cap respectively, see \cite{HL2}. In the following
$\o_d$ denotes the volume of the unit $d$-dimensional sphere in
$\R^{d+1}$.

Given $\th \in (0, \pi)$, we define the spherical cap $\Sig_\th$
in the following way
$$
\Sig_\th = \{x \in S^n : x_{n+1} \geq \cos \th\}.
$$
One can find with elementary computations that the mean curvature of
$\partial \Sig_\th$ with respect to $\Sig_\th$ endowed with the
standard metric $g_0$ is given by
$$
h_{g_0}(\partial \Sig_\th, \Sig_\th)(x') = h_\th := \frac{\cos
\th}{\sin \th}, \qquad \mbox{ for all } x \in \partial \Sig_\th.
$$
We set for brevity $\ov{K} = K(\ov{x})$ (or $K(\ov{x}')$) and $\ov{H}
= H(\ov{x}')$. We want to choose an appropriate $\th$ in such a way
that some solution of the problem
\begin{equation}\label{eq:p-t}
\begin{cases}
-4 \frac{n-1}{n-2} \Delta v + n(n-1)v = \ov{K} v^{\frac{n+2}{n-2}} &
\mbox{ in } \Sig_\th, \\
\frac{2}{n-2} \frac{\partial v}{\partial \nu} + h_\th v = \ov{H} \,
v^{\frac{n}{n - 2}} & \mbox{ in } \partial \Sig_\th.
\end{cases}
\end{equation}
can be chosen to be a constant $v_\theta$. In this way the problem
transforms into
\begin{equation}\label{eq:p-t-c}
\begin{cases}
n(n-1) v_\th = \ov{K} v_\th^{\frac{n+2}{n-2}} &
\mbox{ in } \Sig_\th, \\
h_\th v_\th = \ov{H} \, v_\th^{\frac{n}{n - 2}} & \mbox{ in }
\partial \Sig_\th.
\end{cases}
\end{equation}
From equation (\ref{eq:p-t-c}) it follows that $\th$ must satisfy the relation
\begin{equation}\label{eq:ctk}
\ov{H} \, \sin \th = \left(
\frac{\ov{K}}{n(n-1)}\right)^{\frac{1}{2}} \cos \th,
\end{equation}
and that $v_\th$ solves
\begin{equation}\label{eq:u-c}
\ov{K} \, v_\th^{\frac{4}{n - 2}} = n(n-1).
\end{equation}
As far as the interior blow up is concerned, we look for a constant
function $\hat{v}$ on the whole sphere. Since $\hat{v}$ solves the
equation
\begin{equation}\label{eq:hatv}
-4 \frac{n-1}{n-2} \Delta \hat{v} + n(n-1)\hat{v} = \ov{K}
\hat{v}^{\frac{n+2}{n-2}}, \qquad  \mbox{ in } S^n,
\end{equation}
and is constant, it must also satisfy
\begin{equation}\label{eq:usn}
\ov{K} \, \hat{v}^{\frac{4}{n - 2}} = n(n-1).
\end{equation}

\

\noindent
{\bf Boundary blow up energy}

\

\noindent We now compute the energy of a boundary blow up. Let
$J_{\ov{K}, \ov{H}}$ be the Euler functional as in \eqref{eq:jkh}
corresponding to $K \equiv K(\ov{x}')$ and $H \equiv H(\ov{x}')$; we
have
$$
J_{\ov{K}, \ov{H}}(v_\th) = \frac{1}{2} n(n-1) \int_{\Sig_\th}
v_\th^2 + \frac{(n-1)}{\tan \th} \int_{\partial \Sig_\th} v_\th^2 -
\frac{\ov{K}}{2^*} \int_{\Sig_\th} v_\th^{2^*} - \ov{H} \, (n-2)
\int_{\partial \Sig_\th} v_\th^{2\frac{n-1}{n-2}}.
$$
Taking into account the fact that $v_\th$ is a critical point of
$J_{\ov{K}, \ov{H}}$, and in particular that $J_{\ov{K},
\ov{H}}'(v_\th)[v_\th] = 0$, it turns out that
$$
n(n-1) \int_{\Sig_\th} v_\th^2 + 2 \frac{(n-1)}{\tan \th}
\int_{\partial \Sig_\th} v_\th^2 - \ov{K} \int_{\Sig_\th} v_\th^{2^*}
- 2 \ov{H} (n-1) \int_{\partial \Sig_\th} v_\th^{2\frac{n-1}{n-2}} =
0.
$$
Hence it follows that $$
J_{\ov{K}, \ov{H}}(v_\th) = \left( \frac{1}{2} - \frac{1}{2^*}
\right) \ov{K} \int_{\Sig_\th} v_\th^{2^*} + \ov{H} \int_{\partial
\Sig_\th} v_\th^{2\frac{n-1}{n-2}}.
$$
Setting
$$
F(\th) = \int_0^\th \sin^{n-1} s \, ds; \qquad \th \in [0, \pi],
$$
one immediately checks that
$$
|\Sig_\th| = \o_{n-1} \cdot F(\th); \qquad |\partial \Sig_\th| =
\frac{d}{d\th} |\Sig_\th| = \o_{n-1} \sin^{n-1} \th.
$$
So, taking into account formula \eqref{eq:u-c}, $J_{\ov{K},
\ov{H}}(v_\th)$ can be written as
\begin{equation}\label{eq:ikk}
J_{\ov{K}, \ov{H}}(v_\th) = \frac{\o_{n-1}}{n} \ov{K} F(\th) \cdot
\left( \frac{n(n-1)}{\ov{K}} \right)^{\frac{n}{2}} + \ov{H} \,
\o_{n-1} \sin^{n-1} \th \cdot \left( \frac{n(n-1)}{\ov{K}}
\right)^{\frac{n-1}{2}}.
\end{equation}
Finally, using equation
\eqref{eq:ctk} we deduce
\begin{equation}\label{eq:ebb}
J_{\ov{K}, \ov{H}}(v_\th) = \o_{n-1} \, \left( \frac{n(n-1)}{\ov{K}}
\right)^{\frac{n}{2}-1} \left[ (n-1) \, F(\th) + \cos \th \sin^{n-2}
\th \right].
\end{equation}

\

\noindent
{\bf Interior blow up energy}

\

\noindent Solutions of \eqref{eq:hatv} are critical points of the
functional $J_{\ov{K}} : H^1(S^n) \to \R$ defined by
$$
J_{\ov{K}}(v) = \frac{1}{2} \, \int_{S^n} \left( \frac{4(n-1)}{n-2}
|\n v|^2 + n(n-1) v^2 \right) - \frac{\ov{K}}{2^*}
 \int_{S^n} |v|^{2^*}, \qquad v \in H^1(S^n).
$$
If $\hat{v}$ is the constant given by formula (\ref{eq:usn}), then
its energy is
$$
J_{\ov{K}}(\hat{v}) = \frac{1}{2} n(n-1) \int_{S^n} \hat{v}^2 -
\frac{1}{2^*} \int_{S^n} \ov{K} \hat{v}^{2^*}.
$$
Since $\hat{v}$ is a critical point of $J_{\ov{K}}$ it turns out that
the following relation must be also satisfied
$$
n(n-1) \int_{S^n} \hat{v}^2 -  \int_{S^n} \ov{K} \hat{v}^{2^*} = 0.
$$
Hence it follows that
\begin{equation}\label{eq:ebi}
J_{\ov{K}}(\hat{v}) = \left( \frac{1}{2} - \frac{1}{2^*} \right)
\int_{S^n} \ov{K} \hat{v}^{2^*} = \frac{\o_n}{n} \left( n(n-1)
\right)^{\frac{n}{2}} (\ov{K})^{-\frac{n-2}{2}}.
\end{equation}

\

\noindent
{\bf Comparison of energies}

\

\noindent We conclude this section by proving that, for the same
value of $\ov{K}$, the interior blow up energy is always greater than
the boundary blow up energy, namely we show that
\begin{equation}\label{eq:emi}
J_{\ov{K}, \ov{H}}(v_\th) < J_{\ov{K}}(\hat{v}).
\end{equation}
From equation \eqref{eq:ctk} and from the obvious relation
$$
\o_n = \o_{n-1} \, F(\pi),
$$
one deduces that
\begin{equation}\label{eq:ibue}
J_{\ov{K}}(\hat{v}) = \o_{n-1} \left[ (n-1) \tan^{2-n} \th \, F(\pi)
\right] \, \ov{H}^{2-n}.
\end{equation}
Taking into account \eqref{eq:ctk} and \eqref{eq:ikk}, showing
inequality \eqref{eq:emi} is equivalent to prove
\begin{equation}\label{eq:g-th}
G(\th) := F(\th) + \frac{1}{n-1} \sin^{n-2} \th \cos \th < F(\pi).
\end{equation}
Since it is clearly $G(\pi) = F(\pi)$, we are done if we prove that
$G'(\th) > 0$ for all $\th \in (0, \pi)$. There holds
$$
G'(\th) = \sin^{n-1} \th + \frac{1}{n-1} \left( (n-2)\sin^{n-3} \th \cos^2
\th - \sin^{n-1} \th \right).
$$
With straightforward computations one finally finds that
$$
G'(\th) = \frac{n-2}{n-1} \sin^{n-3} \th \cdot \left( \sin^2 \th +
\cos^2 \th \right),
$$
hence equation (\ref{eq:g-th}) is proved.

We also note that in the case when $\ov{H} = 0$, the boundary blow up
energy is exactly one half of the interior blow up energy, see
formulas \eqref{eq:ebb} and \eqref{eq:ebi}.

\begin{rem}\label{r:phi}
We are particularly interested in the boundary blow up case for $n =
3$. In this situation, using elementary trigonometric formulas, the
explicit expression of $J_{\ov{K}, \ov{H}}(v_\th)$ becomes
$$
J_{\ov{K}, \ov{H}}(v_\th) = 4 \pi \, \sqrt{\frac{6}{\overline{K}}} \,
\left( \frac{\pi}{2} - \arctan\left( \ov{H} \,
\sqrt{\frac{6}{\overline{K}}} \right) \right).
$$
The above function will play a crucial role in the blow up analysis
performed later.
\end{rem}

\section{Blow up analysis: definitions and preliminary results}\label{s:bp}

In this section we recall the definition of isolated and isolated
simple blow up due to R. Schoen, \cite{s}; we also collect some
useful tools and known results.

For a smooth bounded domain $\O \subseteq \R^n$ set $\O_+ = \O
\cap \{ x_n > 0 \}$, $\partial_1 \O = \O \cap \partial \R^n_+$ and
$\partial_2 \O = \partial \O \cap \R^n_+$, hence $\partial \O_+ =
\ov{\partial_1 \O} \cup \partial_2 \O$. We also assume that
$\partial \O$ and $\partial \R^n_+$ intersect transversally, so
that $\partial \O \cap \partial \R^n_+$ is a smooth manifold of
dimension $n - 2$. Let $\nu$ denote the unit exterior normal to
$\O$, and let $\nu'$ denote the exterior unit normal of
$\partial_1 \O$ in $\partial \R^n_+$. Given $w : \partial \R^n_+
\to \R$, the expression $\n' w$ stands for the gradient in
$\R^{n-1}$. If $w$ is defined on $\ov{\O_+}$, the same symbol will
be used for the gradient of the restriction of $w$ to $\partial_1
\O$. In the following $B_\s(x)$ denotes the open ball in $\R^n$ of
radius $\s$ centered at $x$; we just write $B_\s$ if $x = 0$.

We will consider equation \eqref{eq:rp} restricted to $\O$, or
equation \eqref{eq:rpp} when the exponent $p$ is converging to
$\frac{n+2}{n-2}$. For this reason we will not keep the functions $K$
and $H$ fixed, but we will allow them to vary; more precisely, we
consider positive solutions $u_i$ of the sequence of problems
\begin{equation}\label{eq:pb-ap}
\begin{cases}
- \D u_i = \frac{n-2}{4(n-1)} K_i(x) u_i^{p_i}, & \mbox{ in } \O_+; \\
- \frac{\partial u_i}{\partial x_n} = \frac{n-2}{2} H_i(x')
u_i^{\frac{p_i+1}{2}}, & \mbox{ on } \partial_1 \O.
\end{cases}
\end{equation}
We are interested in the case where the supremum of the functions
$u_i$ is tending to infinity, trying to give a precise
characterization of the blow up phenomenon, as in \cite{s} and
\cite{yy1}. A typical ingredient of blow up analysis of scalar
curvature equations is a Pohozahev type identity, which we provide in
the next Lemma.

\begin{lem}\label{l:poh}
Let $p \geq 1$, let $\O \subseteq \R^n$ be as above, and let $K \in
C^1(\ov{\O_+})$, $H \in C^1(\partial_1 \O)$. Assume $u \in
C^2(\ov{\O_+})$ is a positive solution of
\begin{equation}\label{eq:pohe}
\begin{cases}
- \D u = \frac{n-2}{4(n-1)} K(x) \, u^p, \qquad \mbox{ in } \O; \\
- \frac{\partial u}{\partial x_n} = \frac{n-2}{2} H(x') \,
u^{\frac{p+1}{2}}, \qquad \mbox{ on } \partial _1 \O.
\end{cases}
\end{equation}
Then there holds
\begin{eqnarray}\label{eq:poh}
& & \frac{n-2}{4(n-1)} \left(\frac{n-2}{2} - \frac{n}{p+1}\right)
\int_{\O_+} K u^{p+1} + \frac{n-2}{2} \, \left( \frac{n-2}{2} -
\frac{2(n-1)}{p+3} \right) \, \int_{\partial_1 \O} H \,
u^{\frac{p+3}{2}} \nonumber \\ & = & \frac{n-2}{4(n-1)} \frac{1}{p+1}
\, \int_{\O_+} (x \cdot \n K) \, u^{p+1} + \frac{n-2}{p+3} \,
\int_{\partial_1 \O} (x' \cdot \n' H) \, u^{p+1}
\\ & & + \int_{\partial_2 \O} B - \frac{n-2}{4(n-1)} \,
\int_{\partial_2 \O} K \, u^{p+1} \, x \cdot \nu - \frac{n-2}{p+3}
\int_{\partial (\partial_1 \O)} H \, u^{\frac{p+3}{2}} \, x' \cdot
\nu' \nonumber
\end{eqnarray}
where
$$
B = B(x,u,\n u) = \frac{\partial u}{\partial \nu} x
\cdot \n u + \frac{n-2}{2} u
\frac{\partial u}{\partial \nu} - \frac{|\n u|^2}{2} x \cdot \nu.
$$
\end{lem}

\begin{pf}
Multiply the first equation in \eqref{eq:pohe} by $\sum_{j=1}^n \,
x_j \, \frac{\partial u}{\partial x_j}$ and integrate by parts: we
obtain
\begin{eqnarray*}
\frac{n-2}{4(n-1)} \left(\frac{n-2}{2} - \frac{n}{p+1}\right)
\int_{\O_+} K u^{p+1} & - & \frac{n-2}{4(n-1)} \frac{1}{p+1} \,
\int_{\O_+} (x \cdot \n  K) \, u^{p+1} \\ & = & \int_{\partial_2 \O}
B - \frac{n-2}{4(n-1)} \, \int_{\partial_2 \O} K \, u^{p+1} \, x
\cdot \nu. \nonumber
\end{eqnarray*}
Integrating by parts on
$\partial_1 \O$, we deduce
$$
\int_{\partial_1 \O} H \, u^{\frac{p+1}{2}} \, x' \cdot \n'  u = - 2
\, \frac{n-1}{p+1} \, \int_{\partial_1 \O} u^{\frac{p+3}{2}} \, H -
\frac{2}{p+1} \int_{\partial_1 \O} u^{\frac{p+3}{2}} x' \cdot \n' H +
\int_{\partial (\partial_1 \O)} x' \cdot \nu' \, H \,
u^{\frac{p+3}{2}},
$$
Using the second equation in \eqref{eq:pohe}, we easily reach the
conclusion.
\end{pf}

We have also the following Proposition, which proof is elementary.

\begin{pro}\label{c:poh}
Suppose the function $h : (B_\s)_+ \setminus \{0\} \to \R$ is of
the form
$$
h(x) = a \, |x|^{2-n} + b(x),
$$
with $a > 0$ and $b(x)$ of class $C^1$ on $\ov{(B_\s)_+}$. Then there
holds
$$
\lim_{\s \to 0} \int_{\partial_2 B_\s} B(x, h, \n h) = -
\frac{(n-2)^2}{4} \, \o_{n-1} \, a \, b(0).
$$
\end{pro}
Let $\O \subseteq \R^n$ be as above, let $1 < p_i \leq
\frac{n+2}{n-2}$, $p_i \to \frac{n+2}{n-2}$, and let $\t_i =
\frac{n+2}{n-2} - p_i$, so that $\t_i \to 0$. Let $\{K_i\}_i
\subseteq C^1(\ov{\O_+})$, $\{H_i\}_i \subseteq C^1(\ov{\partial_1
\O})$ satisfy for some constant $A_1 > 0$
\begin{equation}\label{eq:A1}
\frac{1}{A_1} \leq K_i(x) \leq A_1, \quad - A_1 \leq H_i(x') \leq
A_1; \qquad \mbox{ for all } x \in \ov{\O_+}, \mbox{ all } x' \in
\ov{\partial_1 \O}, \hbox{ and all } i.
\end{equation}
For every $i$, let also $u_i \in C^2(\ov{\O_+})$ be a positive
solution of problem \eqref{eq:pb-ap}.

\begin{df}\label{d:bu}
The point $\ov{x} \in \O_+ \cup \partial_1 \O$ is called a {\em blow
up point} for $\{u_i\}_i$ if there exists a sequence of points $x_i
\in \O_+ \cup \partial_1 \O$ tending to $\ov{x}$ such that $u_i(x_i)
\to + \infty$.
\end{df}

\begin{df}\label{d:ibu}
Let $\ov{x} \in \O_+ \cup \partial_1 \O$,  and let $\{x_i\}$ be a
sequence of local maxima of $u_i$ such that $x_i \to \ov{x}$ and
$u_i(x_i) \to + \infty$. The point $\ov{x}$ is called an {\em
isolated blow up point} if there exist $0 < \ov{r} < dist(\ov{x},
\partial_2 \O)$ and $\ov{C} > 0$ such that
$$
u_i(x) \leq \ov{C} \, |x - x_i|^{-\frac{2}{p_i - 1}}, \qquad |x -
x_i| \leq \ov{r}, x \in \O_+.
$$
\end{df}
If $\ov{x}$ is a blow up point for $\{u_i\}_i$ we will write for
brevity $x_i \to \ov{x}$ meaning that $\{x_i\}_i$ is a sequence of
points as in Definition \ref{d:ibu}. It is possible to prove,
using Proposition \ref{p:bub} and Lemma \ref{l:harn} below, that
the points $x_i$ having the properties in Definition \ref{d:ibu}
are uniquely determined, provided the functions $K_i$ and $H_i$ in
\eqref{eq:pb-ap} are uniformly bounded in $C^1$ and $C^2$ norm
respectively, see \cite{HL1}.

If $x_i \to \ov{x}$ is a simple blow up for $\{u_i\}_i$ and if
$\ov{r}$ is given by Definition \ref{d:ibu} we define
\begin{equation}\label{eq:uibar}
\ov{u}_i(r) = \frac{1}{|\partial B_r(x_i) \cap \O_+|}
\int_{\partial B_r(x_i) \cap \O_+} u_i, \qquad r \in (0, \ov{r}),
\end{equation}
and
$$
\tilde{u}_i(r) = r^{\frac{2}{p_i - 1}}\ov{u}_i(r), \qquad r \in (0,
\ov{r}).
$$
\begin{df}\label{d:isbu}
The isolated blow up point $x_i \to \ov{x}$ is called {\em isolated
simple} if there exists $\vrh \in (0, \ov{r})$ such that for large
$i$ there holds
\begin{equation}\label{eq:isbu}
\tilde{u}_i \mbox{
has precisely one critical point in } (0,\vrh).
\end{equation}
\end{df}
If $\ov{x}$ is a blow up point, we will call it {\em interior blow up
point} if $\ov{x} \in \O_+$, or {\em boundary blow up point} if
$\ov{x} \in \partial_1 \O$.

\

\noindent Another fundamental tool for the blow up analysis is the
Harnack inequality; we recall the following version from \cite{HL1},
Appendix A.

\begin{lem}\label{l:harn} (Harnack-type inequality)
Let $\{K_i\} \in L^\infty (\O_+)$ and $\{H_i\} \in L^\infty
(\partial_1 \O)$ satisfy \eqref{eq:A1}. Assume also that $\{u_i\}_i$
satisfy (\ref{eq:pb-ap}) with $p_i \geq p_0 > 1$, and let $x_i \to
\ov{x}$ be an isolated blow up point. Then for every $0 < r <
\frac{1}{4} \ov{r}$ the following Harnack-type inequality holds
$$
\sup_{x \in (B_{2r})_+(x_i) \setminus (B_{r/2})_+(x_i)} u_i(x) \leq C
\inf_{x \in (B_{2r})_+(x_i) \setminus (B_{r/2})_+(x_i)} u_i(x),
$$
where $C$ is a positive constant depending only on $n$,  $A_1$ and
$\ov{C}$.
\end{lem}
For the blow up analysis of the first equation in \eqref{eq:pb-ap} we
mainly refer to \cite{yy1}, where the following proposition regarding
the interior blow up points is proved.

\begin{pro}\label{p:old}
Assume $\O \subseteq \R^3$ and that $\{K_i\}_i$ is uniformly
bounded in $C^1(\ov{\O})$. Assume that $p_i \leq \frac{n+2}{n-2}$,
$p_i \to \frac{n+2}{n-2}$, and $\{u_i\}_i$ are solutions of
$$
- \D u_i = \frac{n-2}{4(n-1)} \, K_i \, u_i^{p_i},
\qquad u_i > 0 \hbox{ in } \O.
$$
Then, if $\ov{x} \in \O$ is a blow up point for $u_i$, it is also
an isolated simple blow up point. Moreover there exists an
harmonic function $b : B_{\varrho/2}(\ov{x}) \to \R$ such that,
passing to a subsequence
\begin{equation}\label{eq:ciao}
u_i(x_i) \, u_i(x) \to a \, |x - \ov{x}|^{2-n} + b(x), \qquad
\hbox{ in } C^2_{loc}(B_{\varrho/2} \setminus \{ \ov{x} \}),
\end{equation}
where $a = (4 n (n-1))^{\frac{n-2}{2}}(\lim_i
K_i(x_i))^{\frac{2-n}{2}}$, and where $\varrho$ is given in
Definition \ref{d:isbu}.
\end{pro}

\section{Behavior of isolated simple blow ups}\label{s:bbua}

In this section we perform the study of isolated simple blow ups of
equation \eqref{eq:pb-ap}. The situation of interior blow up has been
treated in \cite{yy1}, hence we are reduced to consider the case in
which the blow up point $\ov{x}$ is in $\partial_1 \O$. We will refer
sometimes to the paper \cite{HL1}, where it is studied equation
\eqref{eq:pb-ap} when $K_i$ and $H_i$ are converging to constant
functions.

\begin{pro}\label{p:bub}
Assume $\{K_i\}_i \subseteq C^1(\ov{\O_+})$ and $\{H_i\}_i \subseteq
C^2(\ov{\partial_1 \O})$ satisfy (\ref{eq:A1}) for some $A_1
> 0$, and satisfy also the condition
\begin{equation}\label{eq:A2}
\| \n K_i\|_{C(\ov{\O_+})} \leq A_2, \qquad \| \n'
H_i\|_{C^1(\ov{\partial_1 \O})} \leq A_2,
\end{equation}
for some $A_2 > 0$. For every $i$, let $u_i$ be a positive solution
of \eqref{eq:pb-ap}, and let $x_i \to \ov{x}' \in
\partial_1 \O$ be an isolated blow up point for $\{u_i\}_i$. Let also
$x'_i$ denote the projection of $x_i$ onto $\partial_1 \O$. Then,
given $R_i \to + \infty$ and $\e_i \to 0^+$, after passing to a
subsequence of $\{u_i\}_i$ (still denoted with $\{u_i\}_i$) we have
$$
\begin{cases}
r_i := u_i(x_i)^{-\frac{p_i-1}{2}} R_i \to 0 \quad \mbox{ as } i
\to + \infty, & \\ \left\|
u_i(x_i)^{-1}u_i(u_i(x_i)^{-\frac{p_i-1}{2}} x + x_i) -
\left(\frac{\l}{(1 + k \l^2(|x'|^2 + |x_n +
t|^2)}\right)^{\frac{n-2}{2}} \right\|_{C^2\left( \left(
u_i(x_i)^{\frac{p_i-1}{2}} (\ov{\O_+} - x_i) \right) \cap \ov{B_{2
R_i}} \right)} \leq \e_i,    &
\end{cases}
$$
where $k = (4n(n-1))^{-1} \, K_i(x'_i)$, while $\l$, $t$ satisfy $2 k
\l t = H_i(x'_i)$ and
\begin{equation}\label{eq:lt}
  \l = \begin{cases}
    1 + k \l^2 t^2 & \text{ if } H_i(x'_i) \geq 0, \\
    1 & \text{ if } H_i(x'_i) < 0.
  \end{cases}
\end{equation}
\end{pro}

\begin{pf}
Consider the functions
$$
w_i(x) = u_i(x_i)^{-1} \, u_i \left( u_i(x_i)^{-\frac{p_i-1}{2}} \, x
+ x_i \right), \qquad x \in u_i(x_i)^{\frac{p_i-1}{2}} (\O_+ - x_i).
$$
It follows immediately that $w_i(0) = 1$ for all $i$ and that $0$ is
a local maximum point for $w_i$. Moreover from the assumption of
isolated blow up we have
$$
w_i (x) \leq C \, |x|^{-\frac{p_i-1}{2}}, \qquad x \in
u_i(x_i)^{\frac{p_i-1}{2}} (\O_+ - x_i) \cap \{ |x| <
u_i(x_i)^{\frac{p_i-1}{2}} \ov{r} \},
$$
where $\ov{r}$ is given in Definition \ref{d:ibu}.

The function $w_i$ is a solution of the problem
$$
\begin{cases}
- \D w_i(x) = \frac{n-2}{4(n-1)} K_i\left(u_i(x_i)^{\frac{p_i-1}{2}}
x + x_i\right) w_i(x)^{p_i}, & \mbox{ in } u_i(x_i)^{\frac{p_i-1}{2}}
(\O_+ - x_i); \\ - \frac{\partial w_i}{\partial x_n}(x) =
\frac{n-2}{2} H_i\left(u_i(x_i)^{\frac{p_i-1}{2}} x + x_i\right) \,
u_i(x)^{\frac{p_i+1}{2}},  & \mbox{ on } u_i(x_i)^{\frac{p_i-1}{2}}
(\partial_1 \O - x_i).
\end{cases}
$$
Denoting by $x_{i,n}$ the $n$-th component of $x_i$ and setting $T_i
= u_i(x_i)^{\frac{p_i-1}{2}} \, x_{i,n}$, two cases may occur, namely
$$
T_i \to + \infty, \qquad \mbox{ or } \qquad T_i \to T \in \R.
$$
In the latter one, we can use \eqref{eq:A1}, \eqref{eq:A2} and the
results in \cite{adn} to prove that the functions $w_i$ converge up
to subsequence, and then one can conclude as in \cite{HL1},
Proposition 1.4.

Hence it is sufficient to exclude the first case. In order to do
this, define the functions
$$
\xi_i(x) = x_{i,n}^{\frac{2}{p_i-1}} \, u_i(x_i + x_{i,n} \, x),
\qquad x \in T_i \, (\O_+ - x_i).
$$
First, letting $\O_i = T_i \, (\O_+ - x_i)$, it is clear that $\O_i$
are relatively open sets which invade the half space $\R^n_1 := \{ x
\in \R^n \, : \, x_n > -1 \}$. Then, since we are supposing by
contradiction that $T_i \to + \infty$, $0$ is an interior blow up
point for the functions $\xi_i$, so from Proposition \ref{p:old} it
follows that $0$ is an isolated simple blow up point. Using Lemma
\ref{l:harn} and the inequality
$$
\xi_i(x) \leq C \, |x|^{-\frac{p_i-1}{2}}, \qquad x \in
T_i^{\frac{p_i-1}{2}} (\O_+ - x_i) \cap \{ |x| <
T_i^{\frac{p_i-1}{2}} \ov{r} \},
$$
the convergence in \eqref{eq:ciao} can be extended to the whole
$\R^n_1 \setminus \{0\}$. Namely one has
$$
\xi_i(0) \, \xi_i(x) \to h(x) \qquad \hbox{ in } C^2_{loc}(\R^n_1
\setminus \{0\}),
$$
where $h(x)$ is a non-negative harmonic function in $\R^n_1 \setminus
\{0\}$ singular at 0 and satisfying
\begin{equation}\label{eq:d1h0}
\frac{\partial h}{\partial x_n} = 0, \qquad
\hbox{ on } \partial \R^n_1.
\end{equation}
By equation \eqref{eq:d1h0} and by the Schwartz's Reflection
Principle, the function $h$ possesses an harmonic extension to the
set $\R^n \setminus \{0, \tilde{0}\}$, where $\tilde{0}$ is the
symmetric point of $0$ with respect to the plane $\partial \R^n_1$.
By uniqueness of harmonic extensions, this must coincide with the
symmetric prolongation of $h$ through $\partial \R^n_1$. Hence the
positivity of $h$ implies that $h(x) = a |x|^{2-n} + A + o(|x|)$ for
$x$ close to $0$ , where $a, A > 0$. Reasoning as in Proposition 3.1
of \cite{yy1}, one can reach a contradiction.
\end{pf}

\noindent Next, we establish the counterpart of Proposition
\ref{p:old} for blow up points in $\partial_1 \O$.

\begin{pro}\label{p:up-barr}
Let $\O = B_2$, suppose $\{K_i\}_i \subseteq C^1(\ov{\O_+})$,
$\{H_i\}_i \subseteq C^2(\ov{\partial_1 \O})$ satisfy conditions
(\ref{eq:A1}) and (\ref{eq:A2}) for some $A_1, A_2 > 0$. Suppose that
for every $i$, $u_i$ satisfies (\ref{eq:pb-ap}) and that $x_i \to 0$
is an isolated simple blow up with
$$
|x - x_i|^{\frac{2}{p_i - 1}}u_i(x) \leq A_3, \qquad \mbox{ for all }
x \in \O_+.
$$
Then there exists some positive constant
$C = C(A_1, A_2, A_3, n, \rh)$ such that
\begin{equation}\label{eq:A4}
u_i(x) \leq C u_i(x_i)^{-1}|x - x_i|^{2-n} \qquad \mbox{ for all } x
\in (B_1(x_i))_+.
\end{equation}
Furthermore, there exists $b : \ov{(B_1)_+}$ satisfying
$$
\begin{cases}
- \D b = 0 & \hbox{ in } (B_1)_+; \\
- \frac{\partial b}{\partial y_n} = 0 & \hbox{ on } B_1 \cap
\partial \R^n_+,
\end{cases}
$$
such that
$$
u_i(x_i)u_i(x) \to a \, |x|^{2-n} + b, \qquad \hbox{ in }
C^2_{loc}\left( \ov{(B_1)_+} \setminus \{0\} \right).
$$
The coefficient $a$ is given by
\begin{equation}\label{eq:aaa}
  a =
  \begin{cases}
    (4n(n-1))^{\frac{n-2}{2}} (\lim_i K_i(x'_i))^{\frac{2-n}{2}}
    & \text{ if } \lim_i H_i(x'_i) < 0;  \\
    \lim_i \left( \frac{K_i(x'_i)}{4n(n-1)} + \frac{H_i(x'_i)^2}{4}
    \right)^{\frac{2-n}{2}} & \text{ if } \lim_i H_i(x'_i) \geq 0,
  \end{cases}
\end{equation}
where $x'_i$ is the projection of $x_i$ onto $\partial_1 \O$.
\end{pro}
Before proving Proposition \ref{p:up-barr}, we need some preliminary
Lemmas.

\begin{lem}\label{l:a-up-barr}
Under the hypotheses of Proposition \ref{p:up-barr}, except for
condition \eqref{eq:A2}, there exist $\d_i > 0$, $\d_i =
O(R_i^{-2+o(1)})$ such that
$$
u_i(x) \leq C \, u_i(x_i)^{-\l_i} \, |x-x_i|^{2-n+\d_i}, \qquad
\mbox{ for } R_i \, u_i(x_i)^{-\frac{p_i-1}{2}} \leq |x-x_i| \leq 1,
$$
where $\l_i = (n - 2 - \d_i)\, \frac{p_i-1}{2} - 1$.
\end{lem}

\begin{pf}
It follows from \cite{HL1}, pages 511-513.
\end{pf}

\begin{lem}\label{l:tia}
Under the hypotheses of Proposition \ref{p:up-barr} there holds
$$
\t_i = O(u_i(x_i)^{-\frac{2}{n-2} + o(1)}), \qquad \hbox{ as } i \to
+ \infty
$$
and therefore
$$
u_i(x_i)^{\t_i} = 1 + o(1), \qquad \hbox{ as } i \to + \infty.
$$
\end{lem}

\begin{pf}
Let $B(x, u, \n u)$ be the function defined in Lemma \ref{l:poh}. By
Lemma \ref{l:a-up-barr}, Proposition \ref{p:bub} and standard
elliptic theories we have
$$
\int_{\partial_2 B_1} B(x,u_i, \n u_i) = O\left( u_i(x_i)^{-2+o(1)}
\right); \qquad \int_{\partial_2 B_1} K_i \, |u_i|^{p_i+1} = O\left(
u_i(x_i)^{-p_i-1+o(1)}\right);
$$
$$
\int_{\partial (\partial_1 B_1)} H_i \, u_i^{\frac{p_i+3}{2}} \, x'
\cdot \nu' = O\left( u_i(x_i)^{-\frac{p_i+3}{2} + o(1)}\right).
$$
Furthermore, since the gradients of $K_i$ and $H_i$ are uniformly
bounded, one can deduce from Lemma \ref{l:a-up-barr}, Proposition
\ref{p:bub} and a  rescaling argument that
$$
\int_{(B_1)_+} u_i^{p_1+1} \, x \cdot \n K_i = O\left(
u_i(x_i)^{-\frac{2}{n-2}+o(1)}\right); \qquad \int_{\partial_1 B_1}
u_i^{\frac{p_i+3}{2}} \, x' \cdot \n' H_i = O\left(
u_i(x_i)^{-\frac{2}{n-2}+o(1)}\right).
$$
On the other hand, using again Lemma \ref{l:a-up-barr} and
Proposition \ref{p:bub} we have also
\begin{eqnarray*}
& & \frac{n-2}{4(n-1)} \left( \frac{n}{p_i+1} - \frac{n-2}{2}
\right) \, \int_{(B_1)_+}  K_i \, u_i^{p_i+1} + \frac{n-2}{2}
\left( \frac{n-2}{2} - 2 \, \frac{n-1}{p_i+3} \right) \,
\int_{\partial_1 B_1} u_i^{\frac{p_i+3}{2}} \, H_i \\ & = & \t_i
\, \frac{(n-2)^3}{16(n-1)} \, \left( \frac{1}{n} K_i(x'_i)
\int_{\R^n_+} \ov{U}_\l^{2^*} + H_i(x'_i) \int_{\partial \R^n_+}
\ov{U}_\l^{2\frac{n-1}{n-2}} + o(1) \right), \qquad i \to +
\infty.
\end{eqnarray*}
Here the function $\ov{U}_\l$ is given by formula \eqref{eq:ulrn+}
with $K(\ov{x}')$ (resp. $H(\ov{x}')$) replaced by $K_i(x'_i)$ (resp.
$H_i(x'_i)$); we note that the values of the above integrals do not
depend on the parameter $\l$.

Using the relation $I'_{K_i(x'_i), H_i(x'_i)}(\ov{U}_\l)[\ov{U}_\l] =
0$, it follows that
$$
\frac{1}{n} K_i(x'_i) \int_{\R^n_+} \ov{U}_\l^{2^*} + H_i(x'_i)
\int_{\partial \R^n_+} \ov{U}_\l^{2\frac{n-1}{n-2}} = I_{\lim_i
K_i(x'_i),\lim_i H_i(x'_i)}(\ov{U}_\l) + o(1) > \a > 0,
$$
where $\a$ is a positive constant depending only on $n$, $A_1$ and
$A_2$. Then the conclusion follows from equation \eqref{eq:poh} and
the above estimates.
\end{pf}

\begin{lem}\label{l:del}
Under the same assumptions of Proposition \ref{p:up-barr} there holds
$$
\lim_i u_i(x_i) \int_{\partial_2 B_1} \frac{\partial u_i}{\partial
\nu} < 0.
$$
\end{lem}

\begin{pf}
Using the divergence formula and \eqref{eq:pb-ap}, we can write
\begin{equation}\label{eq:uin}
u_i(x_i) \int_{\partial_2 B_1} \frac{\partial u_i}{\partial \nu}
= u_i(x_i) \left( \frac{n-2}{4(n-1)} \int_{(B_1)_+} K_i u_i^{p_i} -
\frac{n-2}{2} \int_{\partial_1 B_1} H_i u_i^{\frac{p_i+1}{2}}\right).
\end{equation}
From Lemma \ref{l:a-up-barr} we deduce that
\begin{eqnarray*}
\int_{(B_1)_+ \setminus (B_{r_i}(x_i))_+} u_i^{p_i} & \leq & C \,
\int_{(B_1)_+ \setminus (B_{r_i})_+} \left( u_i(x_i)^{-\l_i} \, |y -
x_i|^{2-n+\d_i} \right)^{p_i} \\ & \leq & C \, R_i^{n-p_i(n-2-\d_i)}
\, u_i(x_i)^{-1+O(\t_i)} = o(1) \, u_i(x_i)^{-1}.
\end{eqnarray*}
In the same way we have that
$$
\int_{\partial _1 B_1 \setminus \partial_1 B_{r_i}(x_i)}
u_i^{\frac{p_i+1}{2}} = o(1) \, u_i(x_i)^{-1}.
$$
Hence, using Proposition \ref{p:bub} and a rescaling argument,
choosing $\e_i \to 0$ sufficiently fast, formula \eqref{eq:uin} can
be written as
$$
u_i(x_i) \int_{\partial_2 B_1} \frac{\partial u_i}{\partial \nu} =
\frac{n-2}{4(n-1)} \lim_i K_i(x'_i) \int_{(B_{R_i})_+}
\ov{U}_{\l}^{\frac{n+2}{n-2}} - \frac{n-2}{2} \lim_i H_i(x'_i)
\int_{\partial_1 B_{R_i}} \ov{U}_{\l}^{\frac{n}{n-2}} + o(1),
$$
where $\l$ satisfies equation \eqref{eq:lt}.
Again from the divergence theorem, we have
$$
\lim_i \left( \frac{n-2}{4(n-1)} K_i(x'_i) \int_{(B_{R_i})_+}
\ov{U}_{\l}^{\frac{n+2}{n-2}} - \frac{n-2}{2}
H_i(x'_i)\int_{\partial_1 B_{R_i}} \ov{U}_{\l}^{\frac{n}{n-2}}
\right) = \lim_{R \to +\infty} \int_{\partial_2 B_R} \frac{\partial
\ov{U}_{\l}}{\partial \nu} < 0.
$$
This concludes the proof.
\end{pf}

\begin{pfn} {\sc of Proposition \ref{p:up-barr}}
Inequality \eqref{eq:A4} for $|x - x_i| \leq r_i$ is an immediate
consequence of Lemma \ref{l:a-up-barr} and Lemma \ref{l:tia}. We now
prove it for $r_i \leq |x - x_i| \leq 1$. Let $\ov{u}_i$ be given by
formula \eqref{eq:uibar} and set $\xi_i(x) = \ov{u}_i(1) \, u_i(x)$.
It is easy to see that $\xi_i$ satisfies
$$
\begin{cases}
- \D \, \xi_i = \frac{n-2}{4(n-1)} \, \ov{u}_i(1)^{p_i-1} \, K_i(x)
\, \xi_i^{p_i}, & \hbox{ in } (B_2)_+; \\
- \frac{\partial \xi_i}{\partial x_n} = \frac{n-2}{2} \,
\ov{u}_i(1)^{\frac{p_i-1}{2}}  H_i(x') \, \xi_i^{\frac{p_i+1}{2}}, &
\hbox{ on } \partial_1 B_2.
\end{cases}
$$
Reasoning as in the proof of Proposition \ref{p:bub} it follows
that, passing to a subsequence, $\{\xi_i\}_i$ converges in
$C^2_{loc}(\ov{(B_2)_+} \setminus 0)$ to some positive function $h
\in C^2_{loc}(\ov{(B_2)_+} \setminus 0)$ satisfying
$$
\begin{cases}
- \D h = 0, & \hbox{ in } (B_2)_+; \\
- \frac{\partial h}{\partial x_n} = 0, & \hbox{ on } \partial_1 B_2.
\end{cases}
$$
Moreover, it follows from condition \eqref{eq:isbu} that $h$ must be
singular at $0$. Reflecting the function $h$ evenly to $B_2$ and
reasoning as above, we deduce that $h$ must be of the form
$$
h(x) = a_1 \, |x|^{2-n} + b(x), \qquad x \in \ov{(B_2)_+} \setminus 0,
$$
where $a_1 > 0$, and $b \in C^2(\ov{(B_2)_+})$ satisfies
$$
\D b(x) = 0, \quad x \in (B_2)_+; \qquad \qquad \frac{\partial
b}{\partial x_n}(x') = 0, \qquad  x' \in \partial_1 B_2.
$$
Now we can prove \eqref{eq:A4} for $|x - x_i| = 1$, namely
\begin{equation}\label{eq:uie}
\ov{u}_i(1) \leq C \, u_i(x_i)^{-1}.
\end{equation}
To do this we observe that, by the harmonicity of $b(y)$ we have
$$
\int_{\partial_2 B_1} \frac{\partial b}{\partial \nu} =
\int_{\partial_1 B_1} \frac{\partial b}{\partial x_n} +
\int_{(B_1)_+} (\D b) = 0,
$$
and so we deduce that
$$
\lim_i \ov{u}_1(1)^{-1} \int_{\partial_2 B_1}
\frac{\partial u_i}{\partial \nu} =
\int_{\partial_2 B_1} \frac{\partial h}{\partial \nu} = a_1
\int_{\partial_2 B_1} \frac{\partial |x|^{2-n}}{\partial \nu} < 0.
$$
Hence formula \eqref{eq:uie} follows from Lemma \ref{l:del}. The
inequality for a general $x$ with $r_i \leq |x - x_i| \leq 1$ follows
from a rescaling argument, as in \cite{yy1} page 340. The value of
the constant $a$ in \eqref{eq:aaa} can be computed multiplying the
first equation in \eqref{eq:pb-ap} by $u_i$, integrating by parts,
and using Proposition \ref{p:bub}.
\end{pfn}

We collect now a couple of technical lemmas which will be needed
later.

\begin{lem}\label{l:s}
Suppose that the hypotheses of Proposition \ref{p:up-barr} hold true.
Then we have the following estimates
$$
\int_{(B_{r_i}(x_i))_+} |x-x_i|^s \, u_i(x)^{p_i+1} = O \left(
u_i(x_i)^{-\frac{2s}{n-2}}\right), \qquad 0 < s < n;
$$
$$
\int_{(B_1(x_i))_+ \setminus (B_{r_i}(x_i))_+} |x-x_i|^s \,
u_i(x)^{p_i+1} = o\left( u_i(x_i)^{-\frac{2s}{n-2}}\right), \qquad
0 < s < n.
$$
$$
\int_{\partial_1 B_{r_i}(x_i)} |x'-x'_i|^s \,
u_i(x')^{\frac{p_i+3}{2}} = O\left(
u_i(x_i)^{-\frac{2s}{n-2}}\right), \qquad 0 < s < n-1;
$$
$$
\int_{\partial_1 B_1(x_i) \setminus \partial_1 B_{r_i}(x_i)}
|x'-x'_i|^s \, u_i(x')^{\frac{p_i+3}{2}} = o \left(
u_i(x'_i)^{-\frac{2s}{n-2}}\right), \qquad 0 < s < n-1;
$$
$$
\int_{\partial_1 B_1(x_i)} |x'-x'_i|^s \, u_i(x')^{p_i+1} =
O\left( u_i(x_i)^{-2\frac{n-1}{n-2}} \, \log u_i(x_i)\right),
\qquad s = n-1.
$$
\end{lem}

\begin{pf}
The proof is a simple consequence of Proposition \ref{p:bub} and
Proposition \ref{p:up-barr}.
\end{pf}

\begin{lem}\label{l:ti}
Suppose that $n = 3$ and that the hypotheses of Proposition
\ref{p:up-barr} hold true. Then we have
$$
\t_i = O(u_i(x_i)^{-2}).
$$
\end{lem}

\begin{pf}
It is sufficient to use \eqref{eq:poh} and Lemma \ref{l:s}.
\end{pf}

\section{Blow up points are isolated simple and at the
boundary}\label{s:6}

In this section we show that blow up points of equation
\eqref{eq:pb-ap} are isolated simple and that the case of interior
blow up can be ruled out. We will use the same terminology about blow
ups for describing both functions on $S^n_+$ or functions defined on
some domain $\O_+ \subseteq \R^n_+$, having in mind the natural
transformation \eqref{eq:iota} induced by the stereographic
projection.

\begin{pro}\label{p:ibisb} Let $n = 3$, let $\O$ be as above,
suppose that $\{K_i\}$ and $\{H_i\}$ satisfy conditions \eqref{eq:A1}
and \eqref{eq:A2}, that for every $i$ $u_i$ is a positive solution of
\eqref{eq:pb-ap} and that $\ov{x}' \in \partial_1 \O$ is an isolated
blow up point for $\{u_i\}$. Then $\ov{x}'$ is also an isolated
simple blow up point.
\end{pro}

\begin{pf}
The proof follows that of Proposition 3.1 in \cite{HL1}, combined
with some argument in \cite{yy1} page 353, and Lemma \ref{l:s}. We
omit the details.
\end{pf}

\begin{pro}\label{p:grt}
Let $\O$, $\{K_i\}$, $\{H_i\}$, $\{u_i\}$ be as in Proposition
\ref{p:ibisb}, and let $\ov{x}'$ be an isolated simple blow up point
for $\{u_i\}$. Suppose also that $\{K_i\}$ are uniformly bounded in
$C^2(\ov{\O_+})$. Let $\varphi_i : \partial_1 B_1 \to \R$ be the
sequence of functions defined by
$$
\varphi_i(x') = 4 \pi \, \sqrt{\frac{6}{K_i(x')}} \, \left(
\frac{\pi}{2} - \arctan\left( H_i(x') \, \sqrt{\frac{6}{K_i(x')}}
\right) \right).
$$
Let also $x'_i$ denote the projection of $x_i$ onto $\partial_1 B_1$.
Then there holds
$$
|\n' \var_i(x'_i) | \leq O(u_i(x_i)^{-2}), \qquad \text{ as } i \to +
\infty .
$$
\end{pro}

\begin{pf}
Choose a test function $\eta \in C^\infty(B_1)$ which satisfies
$$
\eta(x) = 1, \quad x \in B_{1/4}; \qquad \qquad \eta(x) = 0, \quad x
\in B_1 \setminus B_{1/2}.
$$
Multiplying equation \eqref{eq:pb-ap} by $\eta \, \frac{\partial
u_i}{\partial x_j}$, $j = 1, 2$, we obtain
$$
\int_{(B_1)_+} (- \D u_i) \eta \, \frac{\partial u_i}{\partial x_j} =
\frac{1}{8}\int_{(B_1)_+} K_i \, u_i^{p_i} \eta \, \frac{\partial
u_i}{\partial x_j}.
$$
Integrating by parts we deduce
$$
\int_{(B_1)_+} K_i \, u_i^{p_i} \eta \, \frac{\partial u_i}{\partial
x_j} = - \frac{1}{p_i+1} \, \int_{(B_1)_+} u_i^{p_i+1} \left( \eta \,
\frac{\partial K_i}{\partial x_j} + K_i \,  \frac{\partial
\eta}{\partial x_j} \right).
$$
and also
\begin{eqnarray*}
\int_{(B_1)_+} (- \D u_i) \, \eta \, \frac{\partial u_i}{\partial
x_j} & = & \frac{1}{p_i+3} \, \int_{\partial_1 B_1}
u_i^{\frac{p_i+3}{2}} \, \left( \eta \, \frac{\partial H_i}{\partial
x_j} + H_i \, \frac{\partial \eta}{\partial x_j} \right) -
\frac{1}{2} \, \int_{(B_1)_+} |\n u_i|^2 \, \frac{\partial
\eta}{\partial x_j} \\ & & + \int_{(B_1)_+} \n  u_i \, \frac{\partial
u_i}{\partial x_j} \, \n \eta.
\end{eqnarray*}
From the above equations, Proposition \ref{p:up-barr}, and the fact
that $\n \eta$ has support in $\ov{(B_{1/2})_+} \setminus
(B_{1/4})_+$, we obtain
$$
\frac{1}{p_i+1} \, \int_{(B_1)_+} u_i^{p_i+1} \, \eta \,
\frac{\partial K_i}{\partial x_j} + \frac{1}{p_i+3} \,
\int_{\partial_1 B_1} u_i^{\frac{p_i+3}{2}} \, \eta \, \frac{\partial
H_i}{\partial x_j} = O(u_i(x_i)^{-2}).
$$
Using the uniform bounds on the second derivatives of $K_i$ and
$H_i$, and taking into account Lemmas \ref{l:s} and \ref{l:ti} we
deduce
\begin{equation}\label{eq:af2}
\frac{1}{6} \, \frac{\partial K_i}{\partial x_j}(x'_i) \,
\int_{(B_1)_+} u_i^{p_i+1} \, \eta + \frac{1}{8} \, \frac{\partial
H_i}{\partial x_j}(x'_i) \, \int_{\partial_1 B_1}
u_i^{\frac{p_i+3}{2}} \, \eta = O(u_i(x_i)^{-2}).
\end{equation}
Let $\ov{U}_\l$ be the function given in formula \eqref{eq:ulrn+}
with $K(\ov{x}')$ replaced by $K_i(x'_i)$ and $H(\ov{x}')$ replaced
by $H_i(x'_i)$. Using Proposition \ref{p:bub} and Lemma \ref{l:ti},
equation \eqref{eq:af2} becomes
\begin{equation}\label{eq:rel}
\frac{1}{6} \, \frac{\partial K_i}{\partial x_j}(x'_i) \,
\int_{\R^n_+} \ov{U}_\l^{6} + \frac{1}{8} \, \frac{\partial
H_i}{\partial x_j}(x'_i) \, \int_{\partial \R^n_+} \ov{U}_\l^{4} =
O(u_i(x_i)^{-2}).
\end{equation}
By Remark \ref{r:phi}, it turns out that
$$
\var_i(x'_i) = I_{K_i(x'_i),H_i(x'_i)}(\ov{U}_\l) = 4 \,
\int_{\R^n_+} |\n \ov{U}_\l|^2 - \frac{1}{6} \, K_i(x'_i) \,
\int_{\R^n_+} \ov{U}_\l^{6} - H_i(x'_i) \, \int_{\partial \R^n_+}
\ov{U}_\l^{4}.
$$
Differentiating with respect to $x'_i$, and taking into account that
$I'_{K_i(x'_i),H_i(x'_i)}(\ov{U}_\l) = 0$, we deduce that
\begin{eqnarray*}
\frac{\partial \var_i}{\partial y'_j}(x'_i) & = & \frac{\partial
I_{K_i(\cdot),H_i(\cdot)}}{\partial y'_j}|_{x'_i} (\ov{U}_\l) +
I_{K_i(\cdot),H_i(\cdot)} \left( \frac{\partial \ov{U}_{\l(\cdot)}(x,
\cdot)}{\partial y'_j}|_{x'_i} \right) \\ & = & \frac{1}{6} \,
\frac{\partial K_i}{\partial x_j}(x'_i) \, \int_{\R^n_+}
\ov{U}_\l^{6} + \frac{1}{8} \, \frac{\partial H_i}{\partial
x_j}(x'_i) \, \int_{\partial \R^n_+} \ov{U}_\l^{4}.
\end{eqnarray*}
In the above formula the boundary point $y'_i$ is considered as a
parameter on which $I_{K_i,H_i}$ and $\ov{U}_\l$ depend, through the
functions $K_i$ and $H_i$. The conclusion then follows from equation
\eqref{eq:rel} and the last expression.
\end{pf}

\begin{rem}\label{r:f}
We note that if $\{K_i\}$ is just bounded in $C^1$ norm, the above
proof yields anyway $\n' \varphi_i(x'_i) \to 0$ as $i \to + \infty$.
\end{rem}

\noindent Now the local blow up analysis will be applied to equation
\eqref{eq:spa} on the whole half sphere. We begin with the following
Proposition which can be proved as in \cite{HL1} pages 499-502, with
minor modifications.

\begin{pro}\label{p:bu}
Assume $K \in C^1(\ov{S^n_+})$ and $H \in C^{2}(\partial S^n_+)$
satisfy
\begin{equation}\label{eq:kaa}
\frac{1}{A_1} \leq K(x) \leq A_1, \quad \forall x \in S^n_+; \qquad
\qquad \| K \|_{C^1(\ov{S^n_+})} \leq A_2;
\end{equation}
\begin{equation}\label{eq:haa}
- A_1 \leq H(x') \leq A_1, \quad \forall x' \in \partial S^n_+;
\qquad \qquad \| K \|_{C^2(\partial S^n_+)} \leq A_2.
\end{equation}
Then, given any $R \geq 1$ and any $0 < \e < 1$, there exist positive
constants $\d_0, C_0, C_1$, depending only on $n$, $\e, R$, $A_1$ and
$A_2$ such that, for all $\t \leq \d_0$, and for all the solutions
$v$ of equation \eqref{eq:pb-ap} with $\sup_{S^n_+} v \geq C_0$, the
following properties hold true. There exist $\{ q^1, \dots, q^N \}
\subseteq S^n_+$, with $N \geq 1$, such that
\begin{description}
\item{i)} each $q^j$ is a local maximum for $v$ and
$$
\ov{B_{\ov{r}_j}(q^j)} \cap \ov{B_{\ov{r}_l}(q^l)} = \emptyset,
\qquad \hbox{ for } j \neq l,
$$
where $\ov{r}_j = R \, v(q^j)^{-\frac{p-1}{2}}$;
\item{ii)} either $dist(q^j, \partial S^n_+) > \overline{r}_j$ and
\begin{equation*}\tag{$a$}
\left\| v(q^j)^{-1} \, v\left(v(q^j)^{-\frac{p-1}{2}} \, x \right) -
\left( \frac{1}{1 + k_j \, |x|^2} \right)^{\frac{n-2}{2}}
\right\|_{C^2(B_{2R})} < \e,
\end{equation*}
or $dist(q^j, \partial S^n_+) < \overline{r}_j$ and
\begin{equation*}\tag{$b$}
\left\| v(q_i)^{-1} \, v\left( v(q^j)^{-\frac{p-1}{2}} \, x \right) -
\left( \frac{\l_j}{1 + k_j \l_j^2(|x'|^2 + |x_n + t_j|^2)}
\right)^{\frac{n-2}{2}} \right\|_{C^2((B_{2R})_+)} < \e,
\end{equation*}
In the above two formulas it is $k_j = (4n(n-1))^{-1} K(q^j)$, while
$\l_j$ and $t_j$ satisfy $2 k_j \l_j t_j = H(q^j)$, with
$$
\l_j =
  \begin{cases}
    1 + k_j \l_j^2 t_j^2, & \text{ if } H(q^j) \geq 0; \\
    1 & \text{ if } H(q^j) < 0.
  \end{cases}
$$
The function $v$ is identified with its image through the map
$\iota$, the projection being suitably chosen depending on the point
$q^j$.

\item{iii)} $|q^j- q^l|^{\frac{2}{p-1}} \, v(q^l) \geq C_0$,
for $j < l$, while $v(q) \leq C_1 \, \hbox{dist} \left( q, \{ q^1,
\dots, q^N \} \right)^{-\frac{2}{p-1}}$ for all $q \in S^n_+$.
\end{description}
\end{pro}
Properties $(a)$ and $(b)$ in assertion $ii)$ above distinguish,
roughly, the cases of interior and boundary blow ups. Property $iii)$
implies that, if the mutual distance of the points $\{q^j\}$ is
bounded from below along some sequence of solutions, then blow ups
are isolated. This fact is the content of the next Proposition.

\begin{pro}\label{p:dist}
Suppose that $K \in C^1(\ov{S^n_+})$ and $H \in C^2(\partial S^n_+)$
satisfy conditions \eqref{eq:kaa} and \eqref{eq:haa} respectively.
Then, given any $R \geq 1$, and any $0 < \e < 1$, there exist
positive constants $\d_0, \d_1$ and $C_0$ such that, for all $\t \leq
\d_0$, and for all the solutions $v$ of equation \eqref{eq:spa} with
$\sup_{S^n_+} v \geq C_0$ the following property holds true. If $\{
q^1, \dots, q^N \} \subseteq S^n_+$ are the points given by
Proposition \ref{p:bu}, then there holds
$$
\min_{l \neq j} |q^l - q^j| \geq \d_1.
$$
\end{pro}

\begin{pf}
The proof is very similar to that of Proposition 1.2 in \cite{HL1},
and is based on the use of formula \eqref{eq:poh} and on a rescaling
argument. The only difference is that $K$ and $H$ here are non
constant, but one can use conditions \eqref{eq:kaa}, \eqref{eq:haa},
Proposition \ref{p:bub} and Lemma \ref{l:s}.
\end{pf}

\begin{pro}\label{p:noin}
Suppose that $\{K_i\}_i \subseteq C^1(\ov{S^3_+})$ and $\{H_i\}_i
\subseteq C^2(\partial S^3_+)$ satisfy assumptions \eqref{eq:kaa} and
\eqref{eq:haa} uniformly in $i$. Suppose that $\{v_i\}_i$ is a
sequence of positive solutions of \eqref{eq:spa}; then there are no
interior blow-ups for $\{v_i\}_i$.
\end{pro}

\begin{pf}
By Propositions \ref{p:old}, \ref{p:ibisb}, \ref{p:bu}, \ref{p:dist}
we know that both interior and boundary blow ups are isolated simple
and hence isolated. As a consequence, by Definition \ref{d:ibu}, the
number of blow up points is bounded above by a constant depending on
$A_1$ and $A_2$ only.

Suppose by contradiction that $\ov{x} \in S^3_+$ is an interior
blow-up point for $v_i$. Then, it follows from the Harnack
inequality, the fact that there is just a finite number of blow-up
points and Propositions \ref{p:old}, \ref{p:up-barr} that for some
finite set $\{q^1, \dots, q^N\} \subseteq \ov{\R^3_+}$, with $q^1 \in
\R^3_+$, and some harmonic function $b : \R^3_+ \to \R$, the
following holds
$$
u_i(x_i) \, u_i(x) \to a_1 \, |x - q^1|^{-1} + \sum_{j = 2}^N a_j |x
- q^j|^{-1} + b(x), \qquad \hbox{ in } C^2_{loc}(\ov{R}^3_+ \setminus
\{q^1, \dots, q^N\}).
$$
Here $a_j > 0$ for all $j$, $u_i = \iota \, v_i$, and $x_i$ is the
local maximum point of $u_i$ converging to $q^1$ with $u_i(x_i) \to +
\infty$; for simplicity we can suppose that the pole of the
stereographic projection is not a blow up point for $v_i$.

It follows from the Liouville Theorem and from $\frac{\partial
b}{\partial x_3} = 0$ on $\partial \R^3_+$ that $b$ is constant on
$\R^3_+$; reasoning as above we deduce
$$
h(x) = a_0 \, |x - \ov{x}|^{-1} + \ov{b} + O(|x - q^1|), \qquad
\hbox{ for $x$ close to } q^1,
$$
where $\ov{b} > 0$. Let $\s > 0$ be such that $B_\s(q^1) \subseteq
\R^3_+ \setminus \{q^2, \dots, q^N\}$; as for \eqref{eq:poh}, the
function $u_i$ satisfies
\begin{eqnarray}\label{eq:poh1}
\frac{1}{8} \left(\frac{1}{2} - \frac{3}{p_i+1}\right)
\int_{B_\s(q^1)} K_i u_i^{p_i+1} & - & \frac{1}{8(p_i+1)} \,
\int_{B_\s(q^1)} x \cdot \n K_i \, u_i^{p_i+1} + \frac{1}{8} \,
\int_{\partial B_\s(q^1)} K_i \, u_i^{p_i+1} \, x \cdot \nu \nonumber
\\ & = & \int_{\partial B_\s(q^1)} B(x, u_i, \n u_i).
\end{eqnarray}
The estimates of the above terms are completely analogous to the
corresponding ones in boundary blow up analysis, see \cite{yy1}.
Hence, using also Proposition \ref{c:poh}, one deduces that
$$
\left(\frac{1}{2} - \frac{3}{p_i+1}\right) \int_{B_\s(q^1)} K_i
u_i^{p_i+1} = \frac{\t_i}{4} \left( \lim_i K_i(x_i) \int_{\R^3} U_1^6
 + o(1) \right);
$$
$$
\int_{B_\s(q^1)} x \cdot \n K_i(x) \, u_i^{p_i+1} = \int_{B_\s(q^1)}
x \cdot \n K_i(x_i) \, u_i^{p_i+1} + \int_{B_\s(q^1)} x \cdot \n
(K_i(x) - K_i(x_i)) \, u_i^{p_i+1} = o\left( u_i(x_i)^{-2} \right);
$$
$$
\int_{\partial B_\s(q^1)} K_i \, u_i^{p_i+1} = o\left( u_i(x_i)^{-2}
\right); \qquad \int_{\partial B_\s(q^1)} B(x, u_i, \n u_i) = - \pi
\, u_i(x_i)^{-2} a_1 \, \ov{b} + o\left( u_i(x_i)^{-2} \right),
$$
for $\s$ small. Using the last estimates and \eqref{eq:poh1} we reach
a contradiction. This concludes the proof.
\end{pf}

\begin{rem}\label{r:mu}
As anticipated in the Introduction, the fact that there are no
interior blow ups for equation \eqref{eq:spa} is strongly related to
the fact that there are no multiple blow ups for the scalar curvature
equation on $S^3$. The proof relies on the above estimate
$$
\int_{B_\s} x \cdot \n K_i u_i^{\frac{p+1}{2}} =
o\left(u_i(x_i)^{-2}\right).
$$
In our case, the corresponding term $\int_{(B_\s)_+} x \cdot \n
K_i u_i^{\frac{p+1}{2}}$ may be of order $u_i(x_i)^{-2}$ (indeed
this always happens if $(K_i, H_i) \equiv (K, H) \in \mathcal{A}$,
see the proof of Theorem \ref{t:61}) and the above proof breaks
down; this is the reason of the possible presence of multiple blow
ups. See Proposition \ref{t:63}.
\end{rem}

We can summarize the above results with the following proposition.

\begin{pro}\label{p:su}
Suppose that $\{K_i\}_i \subseteq C^1(S^3_+)$ and that $\{H_i\}_i
\subseteq C^2(\partial S^3_+)$ satisfy \eqref{eq:kaa} and
\eqref{eq:haa} uniformly in $i$. Suppose that $\{v_i\}_i$ is a
sequence of positive solutions of
$$
\begin{cases}
- 8 \D v_i + 6 v_i = K_i(x) v_i^{p_i}, & \mbox{ in } S^3_+; \\
2 \frac{\partial v_i}{\partial \nu} = H_i(x') v_i^{\frac{p_i+1}{2}},
& \mbox{ on } \partial S^3_+.
\end{cases}
$$
with $\sup_i v_i = + \infty$. Then the functions $\{v_i\}$ blow up
at a finite number of points of $\partial S^3_+$. These blow ups
are isolated simple and their distance is bounded below by a
positive constant depending on $\min_i K_i$, the $C^1$ bounds of
$\{K_i\}_i$ and the $C^2$ bounds of $\{H_i\}_i$. As a consequence
the number of blow ups is bounded above by a constant depending
only on these numbers. If $(K_i,H_i) = (K,H)$ for some fixed
functions $K$ and $H$, then the blow up points are critical for
$\varphi$.
\end{pro}

%

\section{Proof of Theorem \ref{t:mp}}\label{s:mp}

Consider the subcritical approximation \eqref{eq:spa} of equation
\eqref{eq:sp}, with $p < 5$. From the discussion in Section 2, we are
reduced to find solutions of the equivalent problem \eqref{eq:rpp} in
the half-space. We can choose as pole of the projection the point $-
q$, where $q$ is a global minimum point of $\varphi$, as in the
statement of Theorem \ref{t:mp}. In this way the image of $q$ under
the projection is the origin in $\R^3$.

Solutions of \eqref{eq:rpp} can be found as critical points of the
Euler functional $I_\t : \Dit \to \R$ defined as
$$
I_\t(u) = \frac{1}{2} \int_{\R^3_+} |\n u|^2 - \frac{1}{6 - \t}
\int_{\R^3_+} W^\t(x) K(x) |u|^{6 - \t} - \frac{4}{4 - \t}
\int_{\partial \R^n_+} W^{\frac{\t}{2}}(x') H(x') |u|^{4 - \t},
\qquad u \in \Dit.
$$
Let also $J_\t$ denote the corresponding functional on $H^1(S^3_+)$.
Since the standard half sphere is of {\em positive type} (see the
Introduction), it is clear that the functional $I_\t$ possesses a
mountain pass structure; we denote by $\ov{I}_\t$ the mountain pass
level of $I_\t$. When $\t = 0$, namely when the problem is purely
critical, the functional $I_\t$ is simply $I_{K,H}$, see the notation
in Section \ref{s:p}. We first give an estimate from above of
$\ov{I}_\t$.

\begin{lem}\label{l:emp}
There exist $\d_0$ and $\t_0$, depending on $K$ and $H$, such that
$$
\ov{I}_\t \leq \varphi(q) - \d_0, \qquad \text{ for all } \t \in (0,
\t_0).
$$
\end{lem}

\begin{pf}
For $\l > 0$, let $\ov{U}_\l$ be the function defined in formula
\eqref{eq:ulrn+}, with $k = \frac{K(0)}{4n(n-1)}$ and with $t$
satisfying $2 k \l t = H(0)$. Using a rescaling, it is easy to prove
that
\begin{equation}\label{eq:mp1}
\int_{\R^3_+} |\n \ov{U}_\l|^2 = \int_{\R^3_+} |\n \ov{U}_1|^2;
\qquad \int_{\partial \R^3_+} H(x') \ov{U}_\l^{4} = H(0) \,
\int_{\partial \R^3_+} \ov{U}_1^4 + O\left( \l^{-2} \log \l \right),
\qquad \hbox{ for } \l \hbox{ large},
\end{equation}
and
\begin{equation}\label{eq:mp2}
  \int_{\R^3_+} K(x) \ov{U}_\l^{6} =
  K(0) \int_{\R^3_+} \ov{U}_1^{6} + \l^{-1} \frac{\partial K}{\partial
  x_3}(0)
  \int_{\R^3_+} x_3 \ov{U}_1^{6} + o(\l^{-1}),
  \qquad \hbox{ for } \l \hbox{ large}.
\end{equation}
Using equations \eqref{eq:mp1}, \eqref{eq:mp2} and some simple
computations one finds
$$
\sup_{t \geq 0} I_{K,H}(t \ov{U}_\l) = \varphi(q) - \frac{1}{6}
\l^{-1} \frac{\partial K}{\partial x_3}(0) \int_{\R^3_+} x_3
|\ov{U}_1|^{2^*} + o(\l^{-1}) , \qquad \text{ for } \l \text{ large}.
$$
We note that the condition $\frac{\partial K}{\partial \nu}(q) < 0$
implies $\frac{\partial K}{\partial x_3}(0) > 0$ in $\R^3_+$. Hence,
choosing $\l_0$ sufficiently large, we find the existence of $\d_0$,
depending on $K$ and $H$, such that
$$
\sup_{t \geq 0} I_{K,H}(t \ov{U}_{\l_0}) \leq \varphi(q) - 2 \d_0.
$$
By continuity, choosing $\t_0 > 0$ sufficiently small we deduce that
$$
\ov{I}_\t \leq \sup_{t \geq 0} I_\t(t \ov{U}_{\l_0}) \leq \varphi(q)
- \d_0, \qquad \t \in (0, \t_0).
$$
The proof is thereby completed.
\end{pf}

\

\begin{pfn} {\sc of Theorem \ref{t:mp} concluded.} For $\t > 0$
small, let $v_\t$ be the mountain pass solution of \eqref{eq:spa}. We
claim that the functions $\{v_\t\}_\t$ remain bounded in
$L^\infty(S^3_+)$ as $\t \to 0$. In fact, supposing by contradiction
that $\{v_\t\}_\t$ blows up, by Proposition \ref{p:su} blow ups of
$\{v_\t\}_\t$ occur at the boundary of $S^3_+$ only; let $q^1, \dots,
q^N$ be the blow up points. It follows from Proposition \ref{p:bub}
and Lemma \ref{l:s} that
\begin{equation}\label{eq:suen}
  \lim_{\t \to 0} J_\t(v_\t) = \sum_{j=1}^N \varphi(q^j).
\end{equation}
On the other hand, since $\varphi (q) = \min_{\partial S^3_+}
\varphi$, and since $\varphi > 0$, Lemma \ref{l:emp} contradicts
equation \eqref{eq:suen}.

Hence, the functions $\{v_\t\}$ converge to a solution $v$ of
\eqref{eq:spa}. We note that the function $v$ is non-zero and
strictly positive: this follows from the fact that $v_\t$ is
uniformly away from zero in $H^1(S^3_+)$, or also from the Harnack
inequality. This concludes the proof.
\end{pfn}

\begin{rem}\label{r:n}
Using the Mountain Pass scheme and the standard analysis of Palais
Smale sequences for the functional $I_{K,H}$, one can prove some
existence results of problem \eqref{eq:sp} for any $n$ with more
restrictive hypotheses. Let $\tilde{\varphi}(x')$, $x' \in \partial
S^n_+$, denote the blow up energy corresponding to $K(x')$ and
$H(x')$, computable for example by formula \eqref{eq:ebb} (note that
for $n = 3$ $\tilde{\varphi}$ is nothing but $\varphi$). The
assumptions on $K$ and $H$ are the following. There exists $q \in
\partial S^n_+$ with
\begin{equation}\label{eq:nkb}
\tilde{\varphi}(q) = \min\{ \tilde{\varphi}(x') \, : \, x' \in
\partial S^n_+ \}; \qquad  \frac{\partial K}{\partial \nu}(q) < 0;
\qquad \tilde{\varphi}(q) \leq \frac{\o_n}{n} (n(n-1))^{\frac{n}{2}}
\left( \sup_{S^n_+} K \right)^{-\frac{n-2}{2}} .
\end{equation}
By formula \eqref{eq:ebi}, the last inequality asserts that the
interior blow up has energy larger that the boundary blow up. We also
note that, by \eqref{eq:emi}, this assumption is non empty.

When $n = 3$, the first two conditions in \eqref{eq:nkb} coincide
with \eqref{eq:p-max}, and the last condition can be completely
removed.
\end{rem}

\begin{rem}\label{r:heb}
For the case $H \equiv 0$, Theorem \ref{t:mp} could be proved also
using the observations in Remark \ref{r:n} and the minimization
technique in \cite{hv}.

In fact, if $\sup_{S^3_+} K \leq 4 \, \sup_{\partial S^3_+} K$,
condition \eqref{eq:nkb} is satisfied, see Remark \ref{r:h0}
$(a)$.

On the other hand, if $\sup_{S^3_+} K > \, \sup_{\partial S^3_+}
K$, we can reflect $K$ evenly on all $S^3$ and look for symmetric
solutions of \eqref{eq:scsn}, see the discussion in the
Introduction. Then, using the condition $\sup_{S^3_+} K > 4 \,
\sup_{\partial S^3_+} K$, one can reason as in the proof of
Theorem 4 in \cite{hv}, ruling out concentration of mountain pass
Palais Smale sequences outside $\partial S^3_+$.
\end{rem}

\section{Proof of Theorem \ref{t:61}}\label{s:61}

In this section we prove Theorem \ref{t:61}. We start by giving some
further characterizations of blow up point for solutions of
\eqref{eq:pb-ap}. We recall the definition of the matrix $M_{lj}$
given in formula \eqref{eq:M} and its least eigenvalue $\rho$.

\begin{pro}\label{t:62}
Let $K \in C^1(\ov{S^3_+})$ be a positive function, and let $H \in
C^2(\partial S^3_+)$. Then there exists some number $\d^*
> 0$, depending only on $\min_{\ov{S^3_+}} K$, $\|K\|_{C^1(\ov{S^3_+})}$ and
$\|H\|_{C^2(\partial S^3_+)}$, with the following properties.

Let $\{ p_i \}$ be such that $p_i \leq 5$, $p_i \to 5$, let $K_i \to
K$ in $C^1(\ov{S^3_+})$, $H_i \to H$ in $C^2(\partial S^3_+)$, and
let $v_i
> 0$ satisfy
\begin{equation}\label{eq:paa6}
  \begin{cases}
   - 4\frac{n-1}{n-2} \D v_i = K_i(x) \, v_i^{p_i}  &
   \text{ in } \partial S^3_+, \\
   \frac{2}{n-2} \frac{\partial v_i}{\partial \nu} =
   H_i(x') v_i^{\frac{p_i+1}{2}} & \text{ on } \partial S^3_+.
  \end{cases}
\end{equation}
with $\max_{\ov{S^3_+}} v_i \to + \infty$ as $i \to +\infty$.
Then, after passing to a subsequence, the following properties
hold true
\begin{description}
\item{i)} $\{ v_i \}$ has only isolated simple
blow up points $(q^1, \dots, q^N) \in \mathcal{F} \setminus
\mathcal{F}^-$ $(N \geq 1)$, with $|q^j - q^l| \geq \d^* \; \forall j
\neq l$, and $\rho(q^1, \dots, q^N) \geq 0$. Furthermore $q^1, \dots,
q^N \in  \mathcal{F}^+$ if $N \geq 2$.
\item{ii)} Setting
$$
\mu_j = 2 \, \left( \frac{K(q^j)}{6} + H^+(q^j)^2
\right)^{-\frac{1}{2}} K(q^j)^{\frac{1}{4}} \, \lim_i
\frac{v_i(q_i^1)}{v_i(q_i^l)}
$$
$$
\l_j = \frac{1}{16 \pi} \, \left( \frac{K(q^j)}{6} + H^+(q^j)^2
\right) \frac{\varphi(q^j)}{K(q^j)^{\frac{1}{2}}} \lim_i \t_i \,
v_i(q^j_i)^2,
$$
where $H^+$ is the positive part of $H$ and $q^j_i \to q^j$ is the
local maximum of $v_i$, there holds
$$
\mu_j \in (0, +\infty), \qquad \l_j  \in [0, +\infty), \qquad \forall
j = 1, \dots, N.
$$
\item{iii)} When $N = 1$
$$
\l_1 = \frac{\partial K}{\partial \nu}(q^1) \,
\frac{\psi(q^1)}{K(q^1)^{\frac{3}{2}}};
$$
when $N \geq 2$
\begin{equation}\label{eq:al}
\sum_{l=1}^N M_{lj} \mu_l = \l_j \, \mu_j, \qquad \forall j = 1,
\dots, N.
\end{equation}
\item{iv)} $\l_j \in \, (0,+\infty) \; \forall j = 1, \dots, N$
if and only if $\rho(q^1, \dots, q^N) > 0$.
\end{description}
\end{pro}

\begin{pf}
Assertion $ii)$ follows from Proposition \ref{p:up-barr} and
Lemmas \ref{l:harn}, \ref{l:ti}. From another part, it follows
from Propositions \ref{p:ibisb}, \ref{p:dist} and Remark \ref{r:f}
that $v_i$ has only isolated simple blow up points $q^1, \dots,
q^N \in \mathcal{F}$ $(N \geq 1)$ with $|q^j - q^l| \geq \d^*$ $(j
\neq l)$ for a fixed $\d^* > 0$ depending only on the above
quantities.

Let $q^1_i \to q^1$ be the local maximum of $v_i$ for which
$v_i(q^1_i) \to + \infty$; performing a stereographic projection
through the point $- q^1$, equation \eqref{eq:paa6} is transformed
into
$$
  \begin{cases}
    - \D u_i(x)  = K_i(y) \, W(x)^{\t_i} \, u_i^{p_i} & \text{ in }
    \R^3_+, \\
    - \frac{2}{n-2} \frac{\partial u_i}{\partial x_n} =
    H_i(x') W(x')^{\frac{\t_i}{2}} \, u_i^{\frac{p_i+1}{2}} &
    \text{ on } \partial \R^3_+.
  \end{cases}
$$
By our choice of the projection, it is clear that $0$ is also an
isolated simple blow up point for $\{u_i\}$. We can also suppose that
none of the points $\{q^1, \dots, q^N\}$ is mapped to $+ \infty$ by
the stereographic projection, and we still denote their images by
$q^1, \dots, q^N$ (in particular we have $q^1 = 0$). It follows from
Proposition \ref{p:up-barr} that
\begin{equation}\label{eq:hs6}
u_i(q^1_i) \, u_i(x) \to h_1(x) := a(q^1) \, |x|^{-1} + b_1(x) \qquad
\mbox{ in } C^2_{loc}(\ov{\R^3_+} \setminus \{ q^1, \dots, q^N \}),
\end{equation}
where $a(q^1)$ is the coefficient in \eqref{eq:aaa} with $\lim_i
K_i(x_i)$ replaced by $K(q^1)$ and $\lim_i H_i(x_i)$ replaced by
$K(q^1)$. The function $b_1(x)$ is harmonic in $\ov{\R^3_+}
\setminus \{ q^2, \dots, q^N \}$, and we have still used the
notation $q^j_i$ for the local maxima of $u_i$ converging to
$q^j$.

Coming back to $v_i$ on $S^3_+$ we have
$$
\lim_i v_i(q^1_i) \, v_i(x) = \frac{1}{2} a(q^1) \, G_{q^1}(x) +
\tilde{b}_1(x), \qquad \hbox{ in } C^2_{loc}(\ov{S^3_+} \setminus \{
q^2, \dots, q^N\}),
$$
where $\tilde{b}_1$ is some regular function on $\ov{S^3_+}
\setminus \{ q^2, \dots, q^N \}$ satisfying $(- 8 \D + 6) \,
\tilde{b}_1 = 0$ with $\frac{\partial \tilde{b}_1}{\partial \nu} =
0$ on $\partial S^3_+$.

If $N = 1$, then $\tilde{b}_1 = 0$ by the maximum principle, while
for $N \geq 2$, taking into account the contribution of all the
poles, we deduce that
$$
\lim_i v_i(q^j_i) \, v_i(x) = \frac{1}{\sqrt{2}} a(q^1) \, G_{q^1}(x)
+ \frac{1}{\sqrt{2}} \sum_{l \neq 1} a(q^l) \, G_{q^l}(x) \, \lim_i
\frac{v_i(q^1_i)}{v_i(q^l_i)}.
$$
In fact, subtracting all the poles from the limit function, we
obtain a regular function $r : S^3_+ \to \R$ such that $(- 8 \D +
6) r = 0$ and $\frac{\partial r}{\partial \nu} = 0$ on $\partial
S^3_+$, so it must be $r \equiv 0$. In the above formula, $G_q(x)$
is the function defined in the Introduction, and the convergence
is in $C^2_{loc}(\ov{S^3_+} \setminus \{ q^1, \dots, q^N \})$.

Using the last expression, we can compute the value of $b_1(0)$ in
\eqref{eq:hs6}, which is
\begin{equation}\label{eq:hj}
b_1(0) = \sqrt{2} \, a(q^1) \, \sum_{l \neq 1} a(q^l) \, G_{q^l}(x)
\, \lim_i \frac{v_i(q^1_i)}{v_i(q^l_i)}.
\end{equation}
Hence, using \eqref{eq:hj} and Proposition \ref{c:poh}, we deduce
that
$$
\lim_{\s \to 0}\int_{\partial_2 B_\s} B(x, h_1, \n h_1) = -
\sqrt{2} \, \pi \, a(q^1) \, \sum_{l \neq 1} a(q^l) \,
G_{q^l}(q^1) \, \lim_i \frac{v_i(q^1_i)}{v_i(q^l_i)}.
$$
From another part, it follows from Lemma \ref{l:poh}, Proposition
\ref{p:bub}, Lemma \ref{l:ti} and some computations as in Lemma
\ref{l:del} that
$$
\int_{\partial_2 B_\s} B(x, h_1, \n h_1) = \frac{1}{16}
\varphi(q^1) \, \lim_i \t_i v_i(q_i^1)^2 - \frac{1}{48} \lim_i
u_i(q_i^1)^2 \int_{B_\s} x \cdot \n K_i \, u_i^{p_i+1}.
$$
By Proposition \ref{p:bub}, Lemma \ref{l:ti} and equation
\eqref{eq:intxnu} it follows that
$$
\frac{1}{48} \lim_i u_i(q_i^1)^2 \int_{B_\s} x \cdot \n K_i \,
u_i^{p_i+1} = 3 \pi \, \psi(q^1) \left( 1 + 6
\frac{H^+(q^1)^2}{K(q^1)} \right)^{-1} \, \frac{\partial
K}{\partial x_3}(q^1),
$$
where $\psi$ is defined in \eqref{eq:psi}. The tangent map of the
stereographic projection $\pi$, calculated in $q^1$, is
$\frac{1}{2} \, Id$, hence it turns out that $\frac{\partial (K
\circ \pi^{-1})}{\partial x_3}(q^1) = - 2 \frac{\partial
K}{\partial \nu}(q^1)$. Then, always identifying $K$ with $K \circ
\pi^{-1}$, from the last two formulas we obtain
$$
\frac{1}{16} \varphi(q^1) \, \lim_i \t_i v_i(q_i^1)^2 = - \sqrt{2} \,
\pi \, a(q^1) \, \sum_{l \neq 1} a(q^l) \, G_{q^l}(q^1) \, \lim_i
\frac{v_i(q^1_i)}{v_i(q^l_i)} + \frac{1}{4} \pi \, \frac{\partial
K}{\partial \nu}(q^1) \, a(q^1)^2 \frac{\psi(q^1)}{K(q^1)}.
$$
Finally, using the expression of $\{\mu_l\}$ and $\l_1$ we get
$$
\frac{\psi(q^1)}{K(q^1)^{\frac{3}{2}}} \, \frac{\partial K}{\partial
\nu}(q^1) \, \mu_1 - 4 \sqrt{2} \,
\frac{G_{q^l}(q^1)}{K(q^1)^{\frac{1}{4}} K(q^l)^{\frac{1}{4}}} \mu_l
= \l_1 \, \mu_1.
$$
Of course a similar formula holds for every $q^j$ with $j \neq 1$. We
have thus established \eqref{eq:al} and completed the proof of
$iii)$.

From the last formula it follows that $q^j \in \mathcal{F}
\setminus \mathcal{F}^-$, $\forall j = 1, \dots, N$, and when $N
\geq 2$, $q^j \in \mathcal{F}^+$. Furthermore, since $M_{jj} \geq
0$ for every $j$, and $M_{lj} < 0$ for $l \neq j$, it follows from
linear algebra and the variational characterization of the least
eigenvalue that there exists some $y = (y_1, \dots, y_N) \neq 0$,
$y_l \geq 0 \; \forall l$, such that $\sum_{j = 1}^N \, M_{lj} \,
y_j = \rho \, y_l$.

Multiplying \eqref{eq:al} by $y_j$ and summing over $j$, we have
$$
\rho \, \sum_l \, \mu_l \, y_l = \sum_{l, j} \, M_{lj} \, y_j \,
\mu_l = \sum_j \, \l_j \, \mu_j \, y_j \geq 0.
$$
It follows that $\rho \geq 0$, so we have verified part $i)$. Part
$iv)$ follows from $i)$-$iii)$.
\end{pf}

\

\noindent We introduce now some useful notation. For $\ov{x}' \in
\partial S^3_+$ and $\g > 0$ large, let $\d_{\ov{x}',\gamma} :
S^3_+ \to \R$ be the function defined in the following way
$$
\d_{\ov{x}',\gamma}(x) = \left( \frac{24}{K(\ov{x}')}
\right)^{\frac{1}{4}} \left( \frac{\gamma}{\gamma^2 + 1 + (1 -
\gamma^2) \, \cos d(\tilde{x}, x)} \right)^{\frac{1}{2}},
$$
where $\tilde{x} = (\tilde{x}', \tilde{x}_4) \in S^3$ is given by
$$
\frac{\tilde{x}'}{\left| \tilde{x}' \right|} = \ov{x}'; \qquad \qquad
\tilde{x}_4 = - \sqrt{24} \, \frac{\g}{\g^2 - 1}
\frac{H(\ov{x}')}{\sqrt{K(\ov{x}')}},
$$
and where $d(\tilde{x}, x)$ denotes the geodesic distance in $S^3$.

For all the choices of $\ov{x}'$ and $\g$, $\d_{\ov{x}',\gamma}$
satisfies the equation $(- 8 \D + 6) \d_{\ov{x}',\gamma} =
K(\ov{x}') \d_{\ov{x}',\g}^5$ with the boundary condition $2
\frac{\partial \d_{\ov{x}',\gamma}}{\partial \nu} = H(\ov{x}')
\d_{\ov{x}',\g}^3$. The functions $\{ \d_{\ov{x}',\gamma}
\}_{\ov{x}',\gamma}$, restricted to the half sphere $S^3_+$, are
nothing but the pre-images under the map $\iota$ of the family
$\{\ov{U}_\l\}_\l$ defined in \eqref{eq:ulrn+}, or of some of
their translations in $\R^{n-1}$.

For $\t = p - \frac{n+2}{n-2}$, $\t > 0$, let $J_\t$ denote the Euler
functional corresponding to problem \eqref{eq:spa}, namely
$$
J_\t(v) = 4 \int_{S^3_+} \, |\n v|^2 + 3 \int_{S^3_+} v^2 -
\frac{1}{6-\t} \int_{S^3_+} K(x) \, |v|^{6 - \t} - \frac{4}{4 -
\frac{\t}{2}} \int_{\partial S^3_+} H(x') \, |v|^{4-\frac{\t}{2}},
\qquad v \in H^1(S^3_+).
$$
Let $q^1, \dots, q^N \in \mathcal{F}^+$ be critical points of
$\varphi$ with $\rho(q^1, \dots, q^N) > 0$. For $\e$ small, define
the set $V_{\e} = V_{\e}(\t, q^1, \dots, q^N) \subseteq H^1(S^3_+)$
as
$$
V_{\e} = \left\{ \sum_{i=1}^N \d_{a_i, \gamma_i} \, : \, (\gamma, a)
\in \R^N_+ \times (\partial S^3_+)^N, \quad |a_i - q^i| < \e, \quad
\e < \t \, \gamma_i < \frac{1}{\e}, \quad i = 1, \dots, N \right\}.
$$
We also define $\mathcal{U}_\e = \mathcal{U}_{\e}(\t, q^1, \dots,
q^N)$ to be the $\e$-tubular neighborhood of $V_\e$, namely
$$
\mathcal{U}_\e = \left\{ v + z \, : \, v \in V_\e, z \in (T_v
V_\e)^\perp, \|z\| < \e \right\},
$$
where $(T_v V_\e)^\perp$ denotes the subspace of $H^1(S^3_+)$
orthogonal to $T_v V_\e$.

For $R > 0$, set
$$
\mathcal{O}_R = \left\{ v \in C^{2,\a}(\ov{S^3_+}) \, | \,
\frac{1}{R} \leq v \leq R, \|v\|_{C^{2,\a}(\ov{S^3_+})} \leq R
\right\}.
$$
Using the last definitions and standard regularity results,
Proposition \ref{t:62} can be reformulated as follows.

\begin{pro}\label{p:or}
Let $(K, H) \in \mathcal{A}$ and let $\a \in ]0,1[$. Then there exist
a small positive constant $\e$, and a large positive constant $R$
such that, when $\t > 0$ is sufficiently small, there holds
$$
v \in \mathcal{O}_R \cup \left\{ \mathcal{U}_\e(\t, q^1, \dots, q^N)
\, : \, q^1, \dots, q^N \in \mathcal{F}^+, \rho(q^{1}, \dots q^{N}) >
0, N \geq 1 \right\}
$$
for all $v \in H^1(S^3_+)$ satisfying $v \geq 0$ a.e. and $J'_\t(v) =
0$.
\end{pro}

\noindent Using blow up analysis, we gave necessary conditions on
blowing up solutions of \eqref{eq:spa} when $p$ tends to
$\frac{n+2}{n-2}$ from below. Now we are going to show that if
$(K,H) \in \mathcal{A}$, one can construct solutions highly
concentrating at any $N$ points $q^1, \dots, q^N \in
\mathcal{F}^+$ provided $\rho(q^1, \dots, q^N) > 0$, see
Proposition \ref{t:63} below. The main tool is Implicit Function
Theorem.  Since the procedure is well-known, see \cite{yy2},
\cite{sz}, we just give a general idea of the proof omitting some
details.

We begin with the following technical Lemma, which proof is a
consequence of standard estimates, see \cite{bab}.

\begin{lem}\label{l:a}
Let $\e$ and $\tilde{\d}$ be fixed positive numbers, let $a_i \in
\partial S^3_+$, $i = 1, 2$ be such that $d(a_1, a_2) \geq \tilde{\d}$,
and let $\g_i \in (0, +\infty)$ be such that $\e < \t \g_i <
\frac{1}{\e}$, $i = 1, 2$. Then there exist a positive constant
$C$ such that for $\t$ sufficiently small, the following estimates
hold
$$
\| \d_{a_i, \g_i}^{4-\t} \d_{a_j, \g_j}
\|_{L^{\frac{6}{5}}(S^3_+)} \leq C \, \t, \quad i \neq j; \qquad
\qquad  \| \d_{a_i, \g_i}^{5-\t} - \d_{a_i, \g_i}^5
\|_{L^{\frac{6}{5}}(S^3_+)} \leq C \, \t |\log \t|;
$$
$$
\| d(\cdot, a_i) \d_{a_i, \g_i}^5 \|_{L^{\frac{6}{5}}(S^3_+)} \leq
C \, \t; \qquad \qquad \| \d_{a_i, \g_i}^{3-\t} \d_{a_j, \g_j}
\|_{L^{\frac{4}{3}}(\partial S^3_+)} \leq C \, \t, \quad i \neq j;
$$
$$
\| \d_{a_i, \g_i}^{3-\t} - \d_{a_i, \g_i}^3
\|_{L^{\frac{4}{3}}(\partial S^3_+)} \leq C \, \t |\log \t|;
\qquad \qquad \| d(\cdot, a_i) \d_{a_i, \g_i}^3
\|_{L^{\frac{4}{3}}(\partial S^3_+)} \leq C \, \t.
$$
\end{lem}

Following the original arguments in \cite{bab}, \cite{yy2}, and
using the estimates in Lemma \ref{l:a}, one can prove that for
$\t$ sufficiently small
$$
  \|J'_\t(v)\| \leq O(\t \, |\log \t|), \qquad \hbox{ for } \t \hbox{
  small and } v \in V_{\e}.
$$
Moreover, from Proposition 3.2 in \cite{HL2} and standard
computations, it follows that, for $\t$ small, $I''_\t(u)$ is
invertible in $(T_v V_\e)^\perp$, uniformly with respect to $\t$ and
$v \in V_\e$. Hence by the local inversion theorem, see \cite{ab},
there exists $\e > 0$ small (independent of $\t$) with the following
property. For any $v \in V_\e$, there exists a unique $w(v,\t)$ such
that
\begin{equation}\label{eq:w}
w(v,\t) \in (T_v V_\e)^\perp; \qquad \qquad J'_\t(v + w(v, \t)) \in
T_v V_\e.
\end{equation}
Furthermore, the norm of $w(v, \t)$ can be estimated as
\begin{equation}\label{eq:wut}
\|w(v, \t)\| \leq C \|J'_\t(v)\| \leq C' \, \t \, |\log \t|,
\end{equation}
where $C$ and $C'$ are fixed constants. As a consequence of the above
discussion and of some computations, one finds
\begin{eqnarray}\label{eq:2nd}
J_\t(v + w(v,\t)) & = & J_\t(v) + J'_\t(v)[w(v,\t)] +
O\left(\|w(v,\t)\|^2\right) = J_\t(v) + O\left(\t \, |\log \t|
\right)^2 \nonumber
\\ & = & \sum_{i = 1}^N \varphi(a_i) - \frac{1}{6} \sum_{i = 1}^N
\left( \gamma_i^{-\frac{\t}{2}} - 1 \right) K(a_i) \, \int_{S^3_+}
\d_{a_i, \gamma_i}^6 - \sum_{i = 1}^N \left(\gamma_i^{-\frac{\t}{4}}
- 1 \right) H(a_i) \, \int_{\partial S^3_+}
\d_{a_i, \gamma_i}^4 \\
& & + 4 \, \pi \sqrt{6} \sum_{i = 1}^N \frac{\partial K}{\partial
\nu}(a_i) \frac{\psi(a_i)}{K(a_i)^{\frac{3}{2}}} \,
\frac{1}{\gamma_i} - 64 \, \pi \sqrt{3} \sum_{l \neq j}
\frac{1}{\sqrt{\gamma_j}} \frac{1}{\sqrt{\gamma_j}}
\frac{G_{a_j}(a_l)}{K(a_j)^{\frac{1}{4}} K(a_j)^{\frac{1}{4}}}  +
o(\t), \nonumber
\end{eqnarray}
as $\t \to 0$. By means of equation \eqref{eq:w}, the manifold
$$
\tilde{V}_\e = \left\{ v + w(v, \t), \, : \, v \in V_\e \right\}
$$
is a {\em natural constraint} for $J_\t$, namely a point $u$ which is
critical for $J_\t|_{\tilde{V}_\e}$ is also critical for $J_\t$. In
order to find critical points of $J_\t|_{\tilde{V}_\e}$, we
differentiate $J_\t(v + w(v, \t))$ with respect to the parameters
$a_i$, $\gamma_i$. Using standard estimates we obtain
\begin{equation}\label{eq:dja}
\frac{\partial}{\partial a_i} J_\t(v + w(v, \t)) = \frac{\partial
\varphi}{\partial a_i} + o(1), \qquad v \in V_\e, \t \to 0;
\end{equation}
\begin{eqnarray}\label{eq:djg}
\frac{\partial}{\partial \gamma_j} J_\t(v + w(v, \t)) & = &
\frac{1}{12} \, \frac{\t}{\gamma_j} \, K(a_j) \int_{S^3_+}
\d_{a_j,\gamma_j}^6 + \frac{1}{4} \, \frac{\t}{\gamma_j} \, H(a_j)
\int_{\partial S^3_+} \d_{a_j,\gamma_j}^4 \nonumber
\\ & & - 4 \, \pi \sqrt{6} \sum_{i = 1}^N \frac{\partial K}{\partial
\nu}(a_j) \frac{\psi(a_j)}{K(a_j)^{\frac{3}{2}}} \,
\frac{1}{\gamma_j^2} + 32 \, \pi \sqrt{3} \sum_{l \neq j}
\frac{1}{\gamma_j^{\frac{3}{2}}} \frac{1}{\gamma_l^{\frac{1}{2}}}
\frac{G_{a_j}(a_l)}{K(a_j)^{\frac{1}{4}} K(a_j)^{\frac{1}{4}}} \\ & &
+ o(\t^2), \qquad \qquad v \in V_\e, \t \to 0. \nonumber
\end{eqnarray}
Let us point out that the coefficients of $\frac{1}{\gamma_j^2}$ and
of $\frac{1}{\gamma_j^{\frac{3}{2}}}
\frac{1}{\gamma_l^{\frac{1}{2}}}$ in formula \eqref{eq:djg} coincide,
when $\{a_j\} \equiv \{q^j\}$, with a constant multiple of the
numbers $M_{jj}$ and $M_{lj}$ given in \eqref{eq:M}. As a
consequence, since we are assuming that the least eigenvalue $\rho$
of $(M_{lj})$ is positive, the above coefficients form a positive
definite and invertible matrix. Using this condition, equation
\eqref{eq:dja} and the fact that $\varphi$ is Morse, one can prove
that
\begin{equation}\label{eq:ijve}
  deg_{H^1(S^3_+)} \left( J'_\t|_{\tilde{V}_\e}, \tilde{V}_\e, 0 \right) =
  (-1)^{\sum_{j=1}^N (2 - m(\varphi, q^j))}.
\end{equation}
By the invertibility of $J''_\t$ in the normal direction to $V_\e$,
and by the fact that the functions $\d_{a_i, \gamma_i}$ have Morse
index $1$, it follows from \eqref{eq:ijve} that
\begin{equation}\label{eq:ijue}
    deg_{H^1(S^3_+)} \left( J'_\t, \mathcal{U}_\e,
    0 \right) = (-1)^{N + \sum_{j=1}^N (2 - m(\varphi, q^j))}.
\end{equation}
Since the above degree is always different from zero, $J_\t$ has
at least one critical point in $\mathcal{U}_\e$; moreover it is
standard to prove that critical points of $J_\t$ in
$\mathcal{U}_\e$ are non-negative functions when $\t$ is
sufficiently small. From \cite{adn} and \cite{ch} then it follows
that these solutions are also regular and strictly positive.

We collect the above discussion in the following Proposition.

\begin{pro}\label{t:63}
Let $(K, H) \in \mathcal{A}$, and let $\e > 0$ be small enough.
Then, if $q^1, \dots, q^N \in \mathcal{F}^+$ with $\rho (q^1,
\dots, q^N) > 0$, and if $\t > 0$ is sufficiently small, the
functional $J_\t$ possesses a critical point in
$\mathcal{U}_\e(\t, q^1, \dots, q^N )$. Moreover, formula
\eqref{eq:ijue} holds true and all the critical points of $J_\t$
are strictly positive functions on $S^n$.
\end{pro}
We need now the following lemma which will be useful to obtain a
priori estimates for the computation of some degree formula, see
Proposition \ref{p:cto} below.

\begin{lem}\label{l:hom}
Suppose $(K,H) \in \mathcal{A}$. Then there exists an homotopy
$(K_t,H_t) : [0,1] \to \mathcal{A}$ with the following properties.
\begin{description}
  \item[$(j)$] $H_t = t H$ for all $t \in [0,1]$; moreover $4 \pi
  \sqrt{\frac{6}{K_0}} \equiv \varphi$
  and $K_1 = K$.
  \item[$(jj)$] Setting
$$
\varphi_t(x') = 4 \pi \, \sqrt{\frac{6}{K_t(x')}} \, \left(
\frac{\pi}{2} - \arctan\left( H_t(x') \, \sqrt{\frac{6}{K_t(x')}}
\right) \right), \qquad x' \in \partial S^3_+,
$$
then there holds $\varphi_t \equiv \varphi$ for all $t \in [0,1]$.
  \item[$(jjj)$] Let $\mathcal{F}_t$, $\mathcal{F}_t^\pm$, $\rho_t$
  denote the counterparts of $\mathcal{F}$, $\mathcal{F}^\pm$, $\rho$
  corresponding to the functions $(K_t,H_t)$. Then $\mathcal{F}_t \equiv
  \mathcal{F}$ and $\mathcal{F}_t^\pm \equiv
  \mathcal{F}^\pm$ for all $t \in [0,1]$. Moreover there exists a
  positive constant $\ov{C}$, depending only on $\min_{\ov{S^3_+}}K$,
  $\|K\|_{C^1(\ov{S^3_+})}$, $\|H\|_{C^2(\partial S^3_+)}$ and $\min\{
  |\rho(q^1, \dots, q^N)| \, : \, q^1, \dots, q^N \in \mathcal{F},
  N \geq 1 \}$ such that
\begin{equation}\label{eq:ests}
\min_{\ov{S^3_+}}K_t \geq \frac{1}{\ov{C}}, \qquad
\|K\|_{C^2(\ov{S^3_+})} \leq \ov{C}, \qquad \hbox{ for all } t \in
[0,1];
\end{equation}
\begin{equation}\label{eq:esti}
\qquad \min\{ |\rho(q^1, \dots, q^N)| \, : \, q^1, \dots, q^N \in
\mathcal{F}_t , N \geq 2\} \geq \frac{1}{\ov{C}}, \qquad \hbox{
for all } t \in [0,1].
\end{equation}
\end{description}
\end{lem}

\begin{pf}
First we note that, for a fixed valued of $H(x')$ we have
\begin{equation}\label{eq:limf}
  \lim_{K(x') \to + \infty} \varphi(x') = 0; \qquad
  \lim_{K(x') \to 0} \varphi(x') =
  \begin{cases}
    + \infty & \text{ if } H(x') \leq 0, \\
    \frac{4 \pi}{H(x')} & \text{ if } H(x') > 0.
  \end{cases}
\end{equation}
Moreover, using some simple computations, one finds
\begin{equation}\label{eq:derf}
  \frac{\partial \varphi(x')}{\partial K(x')} = -
  4 \pi \frac{\sqrt{6}}{K(x')^{\frac{3}{2}}} \, \left[
  \frac{\pi}{2} - \arctan \left( H(x') \sqrt{\frac{6}{K(x')}}
  \right) - \frac{H(x') \sqrt{6 K(x')}}{K(x') + 6 H(x')^2} \right]
  < 0.
\end{equation}
As a consequence of \eqref{eq:limf}, \eqref{eq:derf} and the
implicit function theorem, one finds a unique positive function
$K_t(x')$, $x' \in \partial S^3_+$, for which
$$
4 \pi \, \sqrt{\frac{6}{K_t(x')}} \, \left( \frac{\pi}{2} -
\arctan\left( t H(x') \, \sqrt{\frac{6}{K_t(x')}} \right) \right)
= \varphi(x'), \qquad \hbox{ for all } x' \in \partial S^3_+.
$$
We point out that, since $\varphi$ is of class $C^1$, also $K_t$
is of class $C^1$ on $\partial S^3_+$. With such a choice of
$K_t$, properties $(j)$ and $(jj)$ are clearly satisfied.

We are now going prove $(jjj)$, finding a suitable extension of
$K_t$ in the interior of $S^3_+$. Note that by $(jj)$,
$\mathcal{F}_t$ coincides with $\mathcal{F}$ for all $t$. For $q^j
\in \mathcal{F}$, choose $\frac{\partial K_t}{\partial \nu}(q^j)$
satisfying
\begin{equation}\label{eq:dkt}
\frac{\partial K_t}{\partial \nu}(q^j) = \frac{K_t}{K}(q^j) \,
\frac{4 \pi - H(q^j) \varphi(q^j)}{4 \pi - t H(q^j) \varphi(q^j)}
\, \frac{\partial K}{\partial \nu}(q^j), \qquad q^j \in
\mathcal{F}.
\end{equation}
Let $(M_t)_{lj}$ be the counterpart of the matrix $M_{lj}$ defined
in \eqref{eq:M} for the functions $(K_t, H_t)$. It is clear from
\eqref{eq:psi} and \eqref{eq:dkt} that
$$
(M_t)_{lj} = \frac{K(q^l)^{\frac{1}{4}}
K(q^j)^{\frac{1}{4}}}{K_t(q^l)^{\frac{1}{4}}
K_t(q^j)^{\frac{1}{4}}} M_{lj}, \qquad q^j \in \mathcal{F}.
$$
As a consequence, from the multi-linearity of the determinant one
deduces that
$$
\det M_t(q^1, \dots, q^N) = \prod_{j = 1}^N
\frac{K_t(q^j)^{\frac{1}{2}}}{K(q^j)^{\frac{1}{2}}} \det M(q^1,
\dots, q^N),
$$
and hence it follows that $\rho_t(q^1, \dots, q^N) \neq 0$
whenever $\rho(q^1, \dots, q^N) \neq 0$. This implies that $(K_t,
H_t) \in \mathcal{A}$ for all $t$, that $\mathcal{F}_t^\pm \equiv
\mathcal{F}^\pm$ for all $t$, and that \eqref{eq:esti} is
satisfied. Then it is easy to extend $K_t$ in the interior of
$S^3_+$ so that also \eqref{eq:ests} holds true. This concludes
the proof.
\end{pf}

\

\noindent Consider the following problem in $S^3_+$
\begin{equation}\label{eq:eq}
  \begin{cases}
    - 8 \D v + 6 v = f_1 & \text{ in } S^3_+; \\
    2 \frac{\partial v}{\partial \nu} = f_2 & \text{ on } \partial S^3_+.
  \end{cases}
\end{equation}
It is standard, see e.g. \cite{adn}, that if $f_1 \in
C^\a(\ov{S^3_+})$ and if $f_2 \in C^{1,\a}(\partial S^3_+)$ for
some $\a \in (0,1)$, then there exists a solution $v \in C^{2,\a}$
of \eqref{eq:eq}. We denote by $\Xi$ the operator which associates
to $(f_1, f_2)$ the solution $v$ of \eqref{eq:eq}, and we extend
the definition of $\Xi$ also to the case of weak solutions of
\eqref{eq:eq}.

\

\noindent When $(K,H) \in \mathcal{A}$ and the number $\t$ is
bounded from below, we have compactness result for positive
solutions of \eqref{eq:spa} and we can compute their total degree.
We recall the above definition of the set $\mathcal{O}_R$.

\begin{pro}\label{p:cto}
Suppose $(K, 0) \in \mathcal{A}$. Let $J_\t$ denote the Euler
functional corresponding to $(K,0)$. Then there exist constants
$\t_0$, $C_0$ and $\d_0$, depending only on $\min_{S^3_+} K$ and
$\|K\|_{C^1(\ov{S^3_+})}$ with the following properties
\begin{description}
\item{i)}
$\left\{ v \in H^1(S^3_+) \, : \, v \geq 0 \; a.e., \,
J'_{\t_0}(v) = 0 \right\} \subseteq \mathcal{O}_{C_0}$;
\item{ii)}
for $C, \d > 0$ set $\mathcal{O}_{C,\d} = \{ u \in H^1(S^3_+) \, :
\, \exists v \in \mathcal{O}_C$ such that $\|u - v\|_{H^1(S^3_+)}
< \d \}$: then $J'_{\t_0} \neq 0$ on $\partial
\mathcal{O}_{C_0,\d_0}$, and
\begin{equation}\label{dh2}
\mbox{deg}_{H^1(S^3_+)} (u - \Xi(K \, |u|^{4-\t_0}u, 0),
\mathcal{O}_{C_0,\d_0}, 0) = -1.
\end{equation}
\end{description}
\end{pro}

\begin{pf}
Let $\tilde{K} : \ov{S^3_+} \to \R$ be the function defined in the
following way, using stereographic coordinates
$$
\tilde{K}(x) = 2 + \frac{|x|^2 - 1}{x_1^2 + x_2^2 + (x_3 + 1)^2},
\qquad x \in \ov{\R^3_+}.
$$
We point out that $\tilde{K}$ is smooth and strictly positive on
$\ov{S^3_+}$ and satisfies
\begin{equation}\label{eq:kw}
  x \cdot \n \tilde{K}(x) \geq 0, \qquad \hbox{ for all } x \in \ov{\R^3_+}.
\end{equation}
As a consequence, by equation \eqref{eq:poh}, there is no solution
of \eqref{eq:sp} with $(K, H) = (\tilde{K}, 0)$.

Consider the following homotopy from $[0,1]$ into
$C^1(\ov{S^3_+})$
$$
s \to \mathcal{K}_s, \quad \mathcal{K}_s = (1-s) \tilde{K} + s K,
\qquad s \in [0,1].
$$
This homotopy connects $\tilde{K}$ to $K$ when the parameter goes
from $0$ to $1$.

Define $J_{s,\t} : H^1(S^3_+) \to \R$ to be the Euler functional
corresponding to \eqref{eq:spa} for $(K,H) = (\mathcal{K}_s, 0)$.
We claim that for $C_0$ sufficiently large and $\t_0$ sufficiently
small there holds
\begin{eqnarray}\label{eq:bdf}
  \left\{ v \in H^1(S^3_+) \, : \, v \geq 0 \; a.e., J'_{s,\t_0}(v)
  = 0 \right\} \subseteq \mathcal{O}_{C_0},
  \qquad \hbox{ for all } s \in [0,1].
\end{eqnarray}
Of course, by the above discussion, all these weak solutions are
of class $C^{2,\a}$ and positive.

Upper bounds in \eqref{eq:bdf} follow from standard blow up
arguments and from the non-existence results for the problems
\begin{equation}\label{eq:liou}
    - \D u = u^p \qquad \text{ in } \R^n, \qquad u > 0,
\end{equation}
and
\begin{equation}\label{eq:liou2}
  \begin{cases}
    - \D u = u^p & \text{ in } \R^n_+, \qquad u > 0;  \\
    \frac{\partial u}{\partial x_n} = 0 & \text{ on } \partial
    \R^n_+,
  \end{cases}
\end{equation}
when $1 < p < \frac{n+2}{n-2}$. The non-existence result for
\eqref{eq:liou} has been proved in \cite{cgs} while that for
\eqref{eq:liou2} is a consequence of the former.

Once upper bounds are achieved, lower bounds follow from the
Harnack inequality, see Lemma \ref{l:harn}. This proves
\eqref{eq:bdf} and hence property $i)$ in the statement.

Using \eqref{eq:bdf}, it is standard to prove that $J'_{s,\t_0}
\neq 0$ on $\partial \mathcal{O}_{C_0,\d_0}$ for $\d_0$
sufficiently small and for all $s \in [0,2]$; this simply follows
by testing $J'_{s,\t_0}$ on the positive parts of the solutions.

\noindent Therefore, by the homotopy property of the degree, we
only need to establish \eqref{dh2} for $K = \tilde{K}$. In this
case the formula follows from Propositions \ref{t:62}, \ref{p:or}
and \ref{t:63}, since there are no solutions of \eqref{eq:sp} with
$(K, H) = (\tilde{K}, 0)$ and since $\tilde{K}|_{\partial S^3_+}$
possesses just one critical point with $\frac{\partial
\tilde{K}}{\partial \nu} > 0$. This concludes the proof.
\end{pf}

\begin{rem}
In the case in which $H \geq 0$ on $\partial S^3_+$, the a priori
estimates in the previous proof could be obtained from the non
existence results for \eqref{eq:liou} and for the problem
\begin{equation}\label{eq:liou3}
  \begin{cases}
    - \D u = u^p & \text{ in } \R^n_+, \qquad u > 0; \\
    \frac{\partial u}{\partial \nu} = a \, u^q & \text{ on } \partial
    \R^n_+,
  \end{cases}
\end{equation}
where $1 < p < \frac{n+2}{n-2}$, $1< q < \frac{n}{n-2}$ and $a
\leq 0$, see \cite{ccfs} and \cite{lz}. As far as our knowledge,
existence or non-existence of solutions is not known for $a > 0$
and general subcritical exponents $p$ and $q$.
\end{rem}

\begin{pfn} {\sc of Theorem \ref{t:61}}
From the Harnack inequality and standard elliptic estimates it is
enough to prove upper bounds for $v$ in \eqref{eq:cc6}. Arguing by
contradiction, by Proposition \ref{t:62} there exist a sequence of
solutions $\{v_i\}$ of \eqref{eq:sp} blowing up at $q^1, \dots,
q^N \in \partial S^3_+$, and these blow ups are isolated simple.
Taking into account that $(K, H) \in \mathcal{A}$ and that $\l_j =
0$ for all $j$ (since $\t_i = 0$ for all $i$), we get a
contradiction from Proposition \ref{t:62} $iv)$. Hence
\eqref{eq:cc6} is proved.

Let $(\mathcal{K}_t, \mathcal{H}_t)$ be the homotopy defined in
Lemma \ref{l:hom}. We point out that, since the above upper bounds
depend only on $\min_{\ov{S^3_+}}K$, $\|K\|_{C^1(\ov{S^3_+})}$,
$\|H\|_{C^2(\partial S^3_+)}$ and $\min\{ |\rho(q^1, \dots, q^N)|
\, : \, q^1, \dots, q^N \in \mathcal{F}, N \geq 1 \}$, thet are
preserved along the homotopy, by \eqref{eq:ests} and
\eqref{eq:esti}. Hence, using Proposition \ref{p:or} and the
homotopy invariance of the Leray-Schauder degree, we have
\begin{equation}\label{eq:pp1}
\mbox{deg}_{C^{2,\a}(\ov{S^3_+})} (u - \Xi(K \, |u|^{4}u, H \,
|u|^{2}u)), \mathcal{O}_R , 0) = \mbox{deg}_{C^{2,\a}(\ov{S^3_+})}
(u - \Xi(K_0 \, |u|^{4-\t}u, 0), \mathcal{O}_R , 0),
\end{equation}
for $\t$ sufficiently small. Let now $J_\t$ denote the Euler
functional corresponding to $(K_0, 0$. By Propositions \ref{p:or}
and \ref{t:63}, for a suitable value of $\e$ and for $\t$ small,
we know that the non-negative solutions of $J'_\t = 0$ are either
in $\mathcal{O}_R$ or in some $\mathcal{U}_\e(\t, q^1, \dots,
q^N)$; viceversa for all $q^1, \dots, q^N \in \mathcal{F}_+$ with
$\rho(q^1, \dots, q^ k)
> 0$, there are (positive) solutions of $J'_\t = 0$ in
$\mathcal{U}_\e$, and degree of $J'_\t$ on $\mathcal{U}_\e$ is given
by \eqref{eq:ijue}.

Let $C_0 > > R$, $\t_0$ and $\d_0$ be given by Proposition
\ref{p:cto}; take also $\d_1 < < \d_0$. By Proposition \ref{t:63},
\eqref{dh2} and by the excision property of the degree, we have
\begin{equation}\label{eq:pp2}
\mbox{deg}_{H^1(S^3_+)} (u - \Xi(K_0 \, |u|^{4-\t_0}u, 0),
\mathcal{O}_{R,\d_1}, 0) = \mbox{Index}(K, H).
\end{equation}
As in the proof of Proposition \ref{p:cto}, one can check that
there are no critical points of $J_{\t_0}$ in
$\ov{\mathcal{O}_{R,\d_1}} \setminus \mathcal{O}_R$, hence Theorem
B.2 of \cite{yy1} applies and yields
\begin{equation}\label{eq:pp3}
\mbox{deg}_{H^1(S^3_+)} (u - \Xi(K_0 \, |u|^{4-\t_0}u, 0),
\mathcal{O}_{R,\d_1}, 0) = \mbox{deg}_{C^{2,\a}(\ov{S^3_+})} (u -
\Xi(K_0 \, |u|^{4-\t_0}u, 0), \mathcal{O}_R , 0).
\end{equation}
Then the conclusion follows from \eqref{eq:pp1}, \eqref{eq:pp2} and
\eqref{eq:pp3}. The proof of Theorem \ref{t:61} is thereby completed.
\end{pfn}


\begin{thebibliography}{99}

\bibitem{adn} Agmon S., Douglis A., Nirenberg L., Estimates near the
boundary for solutions of elliptic partial differential equations
satisfying general boundary conditions I, Comm. Pure Appl. Math. 12
(1959), 623-727.


\bibitem{ab} Ambrosetti A., Badiale
M., Homoclinics: Poincar\'e-Melnikov type results via a variational
approach, Ann. Inst. Henri. Poincar\'e Analyse Non Lin\'eaire 15
(1998), 233-252. Preliminary note on C. R. Acad. Sci. Paris 323,
S\'erie I (1996), 753-758.


\bibitem{agp} Ambrosetti A., Garcia Azorero J., Peral I., Perturbation of
$-\Delta u + u^{\frac{(N+2)}{(N-2)}} = 0$, the Scalar Curvature
Problem in $\mathbb{R}^N$ and related topics, J. Funct. Anal. 165
(1999), 117-149.


\bibitem{alm} Ambrosetti A., Li Y.Y., Malchiodi A., Yamabe and Scalar
Curvature Problems under boundary conditions, preprint S.I.S.S.A.,
ref. 52/2000/M. Preliminary note on C.R.A.S., S\'{e}rie 1 330 (2000),
1013-1018.


\bibitem{aul} Aubin T., Some Nonlinear Problems in Differential Geometry,
Springer-Verlag, Berlin, (1998).


\bibitem{bab} Bahri A., Critical points at infinity in some variational
problems, Research Notes in Mathematics, 182, Longman-Pitman, London,
1989.

\bibitem{bf} Bahri A., An invariant for Yamabe-type flows with
applications to scalar-curvature problems in high dimension. A
celebration of John F. Nash, Jr. Duke Math. J. 81 (1996),
323--466.


\bibitem{bc} Bahri A., Coron J.  M., The Scalar-Curvature problem on the
standard three-dimensional sphere, J. Funct. Anal. 95 (1991) 106-172.

\bibitem{bcch} Ben Ayed M., Chen Y., Chtioui H., Hammami M., On the
prescribed scalar curvature problem on $4$-manifolds, Duke Math. J.
84 (1996), 633-677.


\bibitem{bp} Bianchi G., Pan X.B., Yamabe equations on
half-spaces, Nonlinear Anal. 37 (1999) 161-186.


\bibitem{cl} Chang K. C., Liu J. Q.,
On Nirenberg's problem, Int. J. Math. 4 (1993), 35-58.

\bibitem{[ACGPY]} Chang S. A., Gursky M. J., Yang P.,
The scalar curvature equation on 2- and 3- spheres, Calc. Var. 1
(1993), 205-229.

\bibitem{cxy} Chang S.A., Xu X., Yang P., A perturbation result for
prescribing mean curvature, Math. Ann. 310-3 (1998), 473-496.

\bibitem{[ACPY1]} Chang S. A., Yang P., Prescribing Gaussian
curvature on $S^2$, Acta Math. 159 (1987), 215-259.

\bibitem{[ACPY2]} Chang, S. A., Yang P., Conformal deformation
of metrics on $S^2$, J. Diff. Geom. 27 (1988), 256-296.

\bibitem{[ACPY3]} Chang S. A., Yang P., A perturbation result
in prescribing scalar curvature on $S^n$, Duke Math. J. 64 (1991),
27-69.

\bibitem{ch} Cherrier P., Probl\'emes de Neumann non lin\'eaires sur les
vari\'et\'es Riemanniennes, J. Funct. Anal.  57 (1984), 154-207.

\bibitem{ccfs} Chipot M., Chlebik M., Fila M., Shafir I.,
Existence of positive solutions of a semilinear elliptic equation
in $\R^n_+$ with a nonlinear boundary condition, J. Math. Anal.
Appl. 223 (1998), 429-471.


\bibitem{dmo} Djadli Z., Malchiodi A., Ould Ahmedou M.,
Prescribing mean curvature on the standard ball, in preparation.


\bibitem{ei} Escobar J., Sharp constant in a Sobolev trace
inequality, Indiana Univ. Math. J. 37 (1988) 687-698.

\bibitem{E1} Escobar J., Conformal deformation of a Riemannian metric
to a scalar fla metric with constant mean curvature on the boundary,
Ann. of Math. 136 (1992), 1-50.

\bibitem{E2} Escobar J., Conformal  metrics with prescribed mean
curvature on the boundary, Cal. Var. 4 (1996), 559-592.

\bibitem{ga} Garcia G., On Conformal Metrics on the Euclidean
Ball, Cornell University, Ph. D.Thesis (2000).

\bibitem{cgs} Gidas B., Spruck J., Global and local behavior
of positive solutions of nonlinear elliptic equations, Comm. Pure
Appl. Math. 34 (1981), 525-598.


\bibitem{HL1} Han Z.C., Li Y.Y., The Yamabe problem on manifolds with
boundaries: existence and compactness results, Duke Math. J. 99
(1999), 489-542.

\bibitem{HL2} Han Z.C., Li Y.Y., The existence of conformal metrics with
constant scalar curvature and constant boundary mean curvature, Comm.
Anal. Geom. 8 (2000), 809-869.


\bibitem{hv} Hebey E., Changements de m\'etriques conformes sur la sph\`ere -
Le probl\`eme de Nirenberg, Bull. Sci. Math. 114 (1990), 215-242.

\bibitem{ll} Li P.L., Liu J.Q., Nirenberg's problem on the
two-dimensional hemi-sphere, Int. J. Math. 4 (1993), 927-939.



\bibitem{yy1} Li Y.Y., Prescribing scalar curvature on $S^{n}$ and
related topics Part I, J.  Diff.  Eq.  120 (1995), 319-410.

\bibitem{yy2} Li Y.Y., Prescribing scalar curvature on $S^{n}$ and
related topics Part II, Existence and compactness, Comm. Pure Appl.
Math. 49 (1996), 437-477.

\bibitem{yn} Li Y.Y., The Nirenberg problem in a domain with
boundary, Top.  Meth.  Nonlin.  Anal.  6 (1995), 309-329.

\bibitem{lz} Li Y.Y., Zhang L., Liouville type theorems and Harnack
type inequalities for semilinear elliptic equations, preprint, 2001.


\bibitem{LZ1} Li Y.Y., Zhu, M.J., Uniqueness theorems through the method of
moving spheres, Duke Math. J. 80 (1995), 383-417.



\bibitem{s} Schoen R., On the number of constant scalar curvature
metrics in a conformal class, in "Differential Geometry: A
Symposium in Honor of Manfredo Do Carmo", (H.B. Lawson and f.
Tenenblat Editors), Wiley, New York, (1991), 331-320.

\bibitem{sy} Schoen R., Yau S.T., Conformally flat manifolds, Kleinian
groups, and scalar curvature, Invent. Math.  92 (1988), 47-71.


\bibitem{sz} Schoen R., Zhang D., Prescribed scalar curvature on
the $n$-sphere, Calculus of  Variations and Partial Differential
Equations, 4 (1996), 1-25.



\end{thebibliography}
\end{document}